\documentclass{amsart}

\title[Pfaffian embeddings of (2,\! 3,\! 5)- into (4,\! 7)-geometries]{A characterization of Pfaffian embeddings from (2,\! 3,\! 5)- into flat (4,\! 7)-geometries}
\thanks{This work was supported by the Institute for Basic Science (IBS-R032-D1).}
\address{Center For Complex Geometry, Institute for Basic Science, 55 Expo-ro, Yuseong-gu 34126 Daejeon, South Korea}
\email{mcmillan@ibs.re.kr}
\author{Benjamin B. McMillan}
\date{April 23, 2024}

\usepackage{amssymb}
\usepackage{amsthm}

\usepackage{mathtools}

\usepackage{microtype}
\usepackage[hidelinks]{hyperref}
\hypersetup{
 pdfauthor={Benjamin McMillan},
 pdftitle={},
 pdfkeywords={},
 pdfsubject={},
 pdfcreator={Emacs 30.0.50 (Org mode 9.6.15)},
 pdflang={English}}

\usepackage{tikz-cd}
\usetikzlibrary{decorations.pathmorphing}

\newcommand{\jhom}{H}

\newcommand{\symdex}[2]{{#1 #2}}
\newcommand{\tendex}[2]{{#1}\cdot{#2}}
\newcommand{\tendexSingle}[1]{{#1}}

\newcommand{\thindex}[1]{{#1}}
\newcommand{\omindex}[1]{{#1}}

\newcommand{\sTheta}[2]{\Theta^{\symdex{#1}{#2}}}
\newcommand{\tTheta}[2]{\Theta^{\tendex{#1}{#2}}}
\newcommand{\tThetaSingle}[1]{\Theta^{\tendexSingle{#1}}}

\newcommand{\ftablA}[3]{\pi^{\symdex{#1}{#2}}_{{#3}}}
\newcommand{\ftablB}[3]{\pi^{\tendex{#1}{#2}}_{{#3}}}
\newcommand{\ftablC}[3]{\pi^{\tendex{#1}{#2}}_{{#3}}}
\newcommand{\ftablBSingle}[2]{\pi^{\tendexSingle{#1}}_{{#2}}}
\newcommand{\ftablCSingle}[2]{\pi^{\tendexSingle{#1}}_{{#2}}}

\newcommand{\tablA}[3]{\jhom^{\symdex{#1}{#2}}_{{#3}}}
\newcommand{\tablB}[3]{\jhom^{\tendex{#1}{#2}}_{{#3}}}
\newcommand{\tablC}[3]{\jhom^{\tendex{#1}{#2}}_{{#3}}}
\newcommand{\tablBSingle}[2]{\jhom^{\tendexSingle{#1}}_{{#2}}}
\newcommand{\tablCSingle}[2]{\jhom^{\tendexSingle{#1}}_{{#2}}}

\newcommand{\Hlltltl}{\jhom^{\symdex{1}{1}}_{{1},{1}}}
\newcommand{\Hlltltz}{\jhom^{\symdex{1}{1}}_{{1},{2}}}
\newcommand{\Hlltztz}{\jhom^{\symdex{1}{1}}_{{2},{2}}}
\newcommand{\Hlztltl}{\jhom^{\symdex{1}{2}}_{{1},{1}}}
\newcommand{\Hlztltz}{\jhom^{\symdex{1}{2}}_{{1},{2}}}
\newcommand{\Hzztltl}{\jhom^{\symdex{2}{2}}_{{1},{1}}}
\newcommand{\Hlbtltl}{\jhom^{\tendex{1}{3}}_{{1},{1}}}
\newcommand{\Hlbtltz}{\jhom^{\tendex{1}{3}}_{{1},{2}}}
\newcommand{\Hlbtlwo}{\jhom^{\tendex{1}{3}}_{{1},{0}}}
\newcommand{\Hlbtlwz}{\jhom^{\tendex{1}{3}}_{{1},{2'}}}
\newcommand{\Hlbptltl}{\jhom^{\tendex{1}{3'}}_{{1},{1}}}
\newcommand{\Hlbptltz}{\jhom^{\tendex{1}{3'}}_{{1},{2}}}
\newcommand{\Hlbptztz}{\jhom^{\tendex{1}{3'}}_{{2},{2}}}
\newcommand{\Hlbptlwo}{\jhom^{\tendex{1}{3'}}_{{1},{0}}}
\newcommand{\Hlbptzwo}{\jhom^{\tendex{1}{3'}}_{{2},{0}}}
\newcommand{\Hlbpwowo}{\jhom^{\tendex{1}{3'}}_{{0},{0}}}
\newcommand{\Hzbtltl}{\jhom^{\tendex{2}{3}}_{{1},{1}}}
\newcommand{\Hzbptltl}{\jhom^{\tendex{2}{3'}}_{{1},{1}}}

\newcommand{\frg}{\mathfrak{g}}

\newcommand{\frm}{\mathfrak{m}}
\newcommand{\frn}{\mathfrak{n}}
\renewcommand{\sl}{\mathfrak{sl}}

\newcommand{\V}{\mathcal{V}}

\newcommand{\BJets}{\mathcal{B} J^{1}}
\newcommand{\B}{\mathcal{B}}
\newcommand{\BP}{\B P}
\newcommand{\Levi}{\operatorname{L}}

\theoremstyle{plain}
\newtheorem{theorem}{Theorem}
\newtheorem*{theorem*}{Theorem}

\newtheorem{proposition}[theorem]{Proposition}

\theoremstyle{definition}

\newtheorem{remark}[theorem]{Remark}

\newtheorem*{repeatedtheorem*}{Theorem \repeatedtheoremnumber}
\newenvironment{theoremduplicate}[1]
  {\def\repeatedtheoremnumber{\ref{#1}}\begin{repeatedtheorem*}}
  {\end{repeatedtheorem*}}

\newcommand{\cblock}[1]{}

\DeclareMathOperator{\Z}{\mathbb{Z}}

\DeclareMathOperator{\R}{\mathbb{R}}

\DeclareMathOperator{\GL}{GL}

\DeclareMathOperator{\Sp}{Sp}

\DeclareMathOperator{\Hom}{Hom}

\DeclareMathOperator{\Sym}{Sym}

\DeclareSymbolFont{script}{U}{eus}{m}{n}
\DeclareMathSymbol{\Wedge}{0}{script}{"5E}

\DeclareMathOperator{\gl}{\mathfrak{gl}}

\renewcommand{\sp}{\mathfrak{sp}}

\DeclareMathOperator{\Ad}{Ad}

\DeclareMathOperator{\I}{\mathcal{I}}
\newcommand{\coframeBundle}{\mathcal{F}}

\renewcommand{\d}{\ensuremath{\hspace{2pt} d}}

\renewcommand{\mod}[1]{\hspace{4mm}\left(\mbox{mod } #1\right)}
\newcommand{\w}{{\mathchoice{\,{\scriptstyle\wedge}\,}{{\scriptstyle\wedge}}
{{\scriptscriptstyle\wedge}}{{\scriptscriptstyle\wedge}}}}

\begin{document}

\begin{abstract}
Given two smooth manifolds with tangent subbundle distributions, an embedding is Pfaffian if its differential sends the distribution on the source into the distribution on the target.
In this paper, we consider the question of existence of Pfaffian embeddings in the specific case where the source is a \( (2,3,5) \)-manifold, the target is the \( 7 \)-dimensional space of isotropic \( 2 \)-planes in a \( 6 \)-dimensional symplectic vector space, and the Pfaffian condition is that the derived \( 3 \)-distribution on the source be mapped into the natural \( 4 \)-distribution on the target.
This is one of the simpler non-trivial cases of the general question on existence of Pfaffian embeddings, but already the answer here requires solution of an interesting differential equation.
It turns out that a generic \( (2,3,5) \)-manifold does not embed, the first obstruction being the fact that a Pfaffian embeddable \( (2,3,5) \)-manifold necessarily has a real double root for its Cartan quartic at each point.
We determine a complete characterization of embeddable \( (2,3,5) \)-manifolds in terms of their associated Cartan geometries, which characterization depends on higher order (non-harmonic) curvature as well.
\end{abstract}
\maketitle
\tableofcontents

\addtocontents{toc}{\protect\setcounter{tocdepth}{1}}

\section{Introduction}\label{sec:Introduction}
\subsection{The problem}\label{sec:The problem}
Given two smooth manifolds \( (M, D) \) and \( (N, D') \) with tangent plane distributions, respectively vector subbundles \( D \subset TM \) and \( D' \subset TN \), one can ask whether there exist embeddings of \( M \) into \( N \) whose differential sends \( D \) into \( D' \).
The question can be stated in a dual form, by replacing \( (M, D) \) and \( (N, D') \) with their equivalent Pfaffian systems \( (M, \I) \) and \( (N, \I') \).
Here we recall that a distribution \( D \) is the same data as a {Pfaffian ideal}, the differentially closed ideal \( \I \) in the graded ring \( \Omega^{*}(M) \) generated locally by any finite collection of \( 1 \)-forms spanning the annihilator of \( D \).
In this dual picture, one asks whether there exists an embedding of \( M \) into \( N \) that pulls back \( \I' \) to be contained in \( \I \).
Both perspectives are useful, but the term Pfaffian is less ambiguous than the term distribution, so we elect to call the maps under consideration \emph{Pfaffian embeddings}.

\subsection{Motivations}\label{sec:Motivations}
There are several reasons that one might ask this question.
It is related to understanding the category of partial differential equations, because there is a natural local correspondence between PDE and Pfaffian systems, and the Pfaffian maps are precisely the maps that preserve the structure of solutions.
The question is also related to the equivalence problem for distributions, because a Pfaffian embedding between systems of equal rank and dimension is a local equivalence.
There has been significant work on the distributional equivalence problem and its generalization, the theory of nilpotent geometries initiated by Tanaka \cite{TanakaDifferentialSystemsGraded1970} and continued by Morimoto \cite{MorimotoGeometricStructuresFiltered1993} and others, 
and it is a natural question after their work to ask how different distributions relate to each other.
This is indeed taken up in the recent paper \cite{Doubrov.Machida.eaExtrinsicGeometryLinear2021} of Doubrov, Machida, and Morimoto, where the authors study the equivalence problem for extrinsic geometries via osculating maps between filtered manifolds.
Their paper sets up a very general framework for studying the local invariants of extrinsic geometries, but doesn't directly answer the question of when such maps exist.

\subsection{Discussion---the specific problem}\label{sec:Discussion---the specific problem}
In this paper, we solve one of the simpler non-trivial cases of this question.
We consider for source \( M \) an arbitrary \( (2,3,5) \)-manifold \( (M, D_{2}) \), meaning that \( D_{2} \) is a generic rank \( 2 \) distribution on a \( 5 \)-manifold, and for target \( N \) the isotropic Grassmannian of \( 2 \)-planes that restrict to zero a constant symplectic form in \( \R^{6} \).
The space \( N \) is a homogeneous space, and has a natural choice of invariant \( 4 \)-plane distribution \( D'_{4} \).
We ask when there exist Pfaffian embeddings \( (M, D_{3}) \to (N, D'_{4}) \) sending the derived \( 3 \)-distribution \( D_{3} = [D_{2}, D_{2}] \) on \( M \) into the \( 4 \)-distribution on \( N \).

To describe why this would be a natural choice, we first observe that \( (M, D_{3}) \) and \( (N, D'_{4}) \) are {regular differential systems}, in the sense of \cite{TanakaDifferentialSystemsGraded1970}.
Any regular differential system has an associated bundle of nilpotent graded Lie algebras---its {symbol bundle}---and a Pfaffian embedding determines a map of these symbol bundles.
The symbol mapping provides a sort of linearization of the embedding question, and in many cases one can rule out the possibility of Pfaffian embeddings at the symbol level, simply by showing the symbols are not compatible.
The symbols of the \( M \) and \( N \) considered here are both of constant type---on \( M \) because there is only one isomorphism class of non-degenerate \( (3,5) \) nilpotent graded Lie algebra,
and on \( N \) because the distribution is invariant under the homogeneous structure of \( N \).
It is then an algebraic computation to see that the symbols of \( M \) and \( N \) are compatible, the space of Lie algebra embeddings between them being a positive dimensional manifold.
So there is no obstruction to embeddings at the symbol level, and it turns out that embeddings do exist for some, but not all, choices of \( (M, D_{3}) \).

The discussion of the previous paragraph holds for any \( 5 \)-manifold \( M \) with a bracket generating \( 3 \)-distribution, and it remains to explain how \( (2,3,5) \)-manifolds enter the discussion.
The point is that the first invariant of a \( 3 \)-distribution on \( M \) is its square-root \( 2 \)-distribution \( D_{2} \), defined as the kernel of the Levi-bracket on \( D_{3} \).
Enforce the mild constant rank assumption that this \( 2 \)-distribution is either everywhere Frobenius or nowhere Frobenius.
In the Frobenius case, \( (M, D_{3}) \) is locally isomorphic to the flat model, \( J^{1}(\R, \R^{2}) \) with its contact distribution, and it follows easily from the compatibility of symbols discussed above that the flat model admits Pfaffian embeddings into \( (N, D'_{4}) \).
So, with little loss of generality, suppose that \( D_{2} \) is nowhere Frobenius, in which case its derived \( 3 \)-distribution equals \( D_{3} \), and so \( D_{2} \) generates a \( (2,3,5) \)-structure.

Once this additional geometric structure is available, one should take advantage of it. 
Indeed, a further reason to consider \( (2,3,5) \)-manifolds as a source is that their equivalence problem is well understood, allowing the possibility of a reasonable statement regarding which ones embed.
It is well known that \( (2,3,5) \)-manifolds are in correspondence with certain regular normal Cartan geometries (see for example \cite{Cap.SlovakParabolicGeometries2009}, or \cite{TheCartantheoreticClassificationMultiplytransitive2022} for a discussion specific to \( (2,3,5) \)), and one can hope to characterize the Pfaffian embeddable manifolds in terms of the curvature of the Cartan geometry.
This turns out to be the case, and we find that \( M \) only admits embeddings if the Cartan geometry admits a reduction satisfying specific conditions.

More precisely, the Cartan geometry over \( M \) is a canonically coframed principal bundle \( \B_{M} \), with structure group \( P \) a certain subgroup of the exceptional simple Lie group \( G_{2} \).
If \( M \) is Pfaffian embeddable, then one can fix the principal subbundle \( \mathcal{R} \) of \( \B_{M} \) comprising coframings adapted to a given embedding.
Such an adapted reduction will be quite special geometrically, and its existence strongly constrains the Cartan curvature of \( M \).
Existence can be characterized via a collection of forms given in terms of the Cartan coframing and the curvature functions of \( \B_{M} \); 
in notation to be explained below (Subsection \ref{sec:Curved (2,3,5)-manifolds}), let
\begin{multline*}
 I_{M} = \{ \gamma^{1}_{2},\enspace
      \gamma^{0}_{2} - \tfrac{1}{14} A_{3} \theta^{1},\enspace
      \gamma + \tfrac{4}{7} B_{3} \theta^{1} + \tfrac{5}{7} A_{3} \omega^{0},\enspace
      \gamma_{2} + \tfrac{17}{7} C_{2} \theta^{1} - \tfrac{17}{14} A_{3} \omega^{1'},\enspace
      \\ \gamma_{1} + \tfrac{5}{7} C_{3} \theta^{1} - C_{2} \theta^{2} + \tfrac{22}{7} B_{3} \omega^{0} - \tfrac{9}{7} A_{4} \omega^{1'} - \tfrac{37}{14} A_{3} \omega^{2'} \} .
\end{multline*}
For \( \mathcal{R} \) to be an adapted reduction it is necessary, and almost sufficient, that these forms restrict to zero on \( \mathcal{R} \).
Here is roughly the main Theorem of the paper, with the full statement given in Subsection \ref{sec:The main theorem}.
\begin{theoremduplicate}{thm: main embeddability theorem}
  The manifold \( (M, D) \) admits Pfaffian embeddings into \( (N, D') \) if and only if \( \B_{M} \) admits a principal reduction \( \mathcal{R} \subset \B_{M} \) such that the forms of \( I_{M} \) restrict to zero and two further conditions hold on the curvature functions restricted to \( \mathcal{R} \).
\end{theoremduplicate}
The existence of an adapted reduction \( \mathcal{R} \) constrains the curvature functions because it determines a privileged direction in the distribution \( D_{2} \), and it is not difficult to see that the Cartan quartic must have a \emph{real} double root in this direction (in Cartan's notation, \( A_{1} = A_{2} = 0 \)), as well as a double root for the Cartan cubic (\( B_{1} = B_{2} = 0 \)), a root for the quadratic (\( C_{1} = 0 \)), a root for the linear part (\( D_{1} = 0 \)), as well as further curvature conditions.
The requirement of a real double root is already strong, implying for example constraints on the Weyl curvature of Nurowski's conformal structure \cite{Nurowski2005--DifferentialEquationsConformalStructures}.

When embeddings do exist, they are given in an explicit form.
In the course of the proof, one arrives at a differential ideal \( \I^{(2)} \) defined on a bundle \( \V \) over \( M \times N \), and the conditions given in Theorem \ref{thm: main embeddability theorem} are precisely that \( \I^{(2)} \) is a Frobenius ideal.
In case \( \I^{(2)} \) is Frobenius, its integral leaves each project down to the graph of a Pfaffian embedding.
The leaves themselves can be regarded as the bundle over this graph of coframes of \( M \) and \( N \) that are adapted to the embedding.

The reduction of the Cartan \( G_{2} \)-structure in Theorem \ref{thm: main embeddability theorem} reflects the extrinsic geometry of the embedding, which is why it is a bundle of fiber dimension \( 4 \) even in case of an embedded flat \( (2,3,5) \).
If the embedded \( (2,3,5) \)-system is not flat then there are further reductions.

In the final Section \ref{sec:Examples} we discuss some examples that embed.
Of particular interest may be the case of a \( (2,3,5) \)-system whose Cartan quartic has two distinct double roots.
It follows from Theorem \ref{thm: main embeddability theorem} and the consequent compatibility conditions that there is a unique such system that embeds.
A posteriori, this embeddable system is homogeneous, and comparison with known classification results \cite{Willse2019--HomogeneousRealDistributionsIsotropy} shows that it is the system obtained from rolling two de Sitter spaces of curvature ratio \( 4:1 \).

\subsection{Discussion---the method}\label{sec:Discussion---the method}
The approach taken to prove Theorem \ref{thm: main embeddability theorem} is simply to solve the partial differential equation (really, exterior differential system) whose solutions correspond to Pfaffian embeddings of \( M \) into \( N \).
The EDS to be solved is the obvious one, defined on the space \( P \) of \( 1 \)-jets of local embeddings \( M \) to \( N \) that send \( D \) to \( D' \).
The contact form \( \Theta \) of \( J^{1}(M, N) \) pulls back to \( P \), and Pfaffian embeddings are in bijection with sections of \( P \to M \) that pull \( \Theta \) back to \( 0 \).
So one can analyze the linear Pfaffian system \( (P, \I) \), where \( \I \) is the differentially closed ideal generated by the component \( 1 \)-forms of \( \Theta \).

In the solution of the EDS, some care is taken to incorporate the well understood geometry of both the source and target, and the symmetries therein.
To this end we replace the usual jet bundle \( J^{1} = J^{1}(M,N) \) with the \emph{coframed jet bundle} \( \BJets = \BJets(M, N) \), which accounts for these symmetries.
The incorporation of coframings into \( \BJets \) trivializes, and thus simplifies the description of, the various bundle-valued objects on \( J^{1} \).
For example, the contact form \( \Theta \) is a \( TN \)-valued \( 1 \)-form on \( J^{1} \), but its pullback to \( \BJets \) takes values in a single vector space.
One can thus fix a basis on this vector space and obtain a canonical choice of components of \( \Theta \).
The structure equations for these components are simple to write explicitly on \( \BJets \), and one has a good choice of tableau forms for the linear Pfaffian system.
In the end, one obtains a canonical (partial) linear Pfaffian coframing on \( \BJets \).

The EDS \( (P, \I) \) can be lifted to a larger system \( (\BP, \I') \), with \( \BP \) a submanifold of \( \BJets \).
This lifted system \( \BP \) is essentially equivalent to the original one, being simply a bundle over \( P \) whose fibers are the Cauchy characteristic leaves of \( \I' \).
In particular, the integral sections of \( P \) correspond to integral submanifolds of \( \BP \), in the usual manner.

This replacement allows one to make uniform reductions of coframe relative the jets of solutions.
For a more familiar analogue of this approach, consider the problem of isometric embeddings into Euclidean space, as done for example in \cite{Bryant.Chern.eaExteriorDifferentialSystems1991}, section VII.4.
There too, the authors restate the question as an EDS on a bundle of coframed jets, and continue the analysis from there.
In that case, it is clear that one should immediately normalize framings on the source and target so that each \( 1 \)-jet looks like a fixed canonical isometry.
It is an important point that these normalizations are made at the level of \( 1 \)-jets of solutions, so that they can be done uniformly, and before actually solving the PDE.

Similar questions of normalization arise here, but the normalizations are more complicated to uncover.
Here a further advantage of working on \( \BJets \) is seen.
While solving the EDS, one must restrict to submanifolds cut out by integrability conditions.
The first such restriction requires one to solve a non-linear system of equations on \( \BP \), but by the trivializations described above it suffices to solve the equations on a single vector space.
Then, at cost of shrinking the Cauchy leaves of \( \BP \), the symmetry groups of the coframings on \( M \) and \( N \) can be used to normalize the solutions to a standard form.
The standard form can be chosen to simplify the new symbol relations that are necessarily introduced by restriction to the integrability condition.
Such reductions are made throughout the calculation, which would be significantly more complicated otherwise.

The approach used here can certainly be applied to the question of Pfaffian embedding for other distributions, or more generally embeddings between nilpotent geometries.
Indeed, the calculation here becomes complicated, but is effective, perhaps with the aid of a good computer algebra package.

We note that it is not essential to the calculation that the symbols of \( M \) and \( N \) are of constant type, and the method used here can be used to find conditions for embeddings between distributions of non-constant symbol type.
For example, in \cite{Hwang.LiLinesHolomorphicContact2023}, Hwang and Li introduce a class of distributions that generalize \( (2,3,5) \)-manifolds in a precise sense, and in particular have small growth vector \( (2m, 3m, 3m+2) \) for \( m \ge 1 \).
For \( m > 2 \), a generic \( (M, D) \) in their class has continuously varying moduli of symbol type.
Nonetheless, the symbol at each point of \( M \) is compatible with the symbol of the natural \( 4m \)-distribution on the homogeneous space \( N \) of isotropic \( 2 \)-planes in a symplectic vector space \( \R^{2m+4} \).
So, again there is no obstruction to embeddings at the symbol level.
The approach of this paper can be used to find invariants on such \( M \) that obstruct Pfaffian embeddings into \( N \).
In fact, the first invariants one finds involve the variation of symbol, vanishing when \( M \) has constant symbol type.
Of course, there may be further obstructions even for \( M \) of constant symbol, as seen in the example given in this paper.

On the other hand, it would be interesting to understand the specific case studied here from a more algebraic or geometric perspective.
Alternative approaches that rely more heavily on constant symbol type may well be more efficient for determining the conditions for embeddable \( (2,3,5) \)-manifolds.
For example, the extrinsic geometric invariants of \cite{Doubrov.Machida.eaExtrinsicGeometryLinear2021} allow one to contemplate Gauss-Codazzi type relations between the extrinsic and intrinsic curvatures, and such a deep understanding of this case would likely recover the conditions found here, with the advantage of a deeper geometric explanation.

\subsection{The structure of the paper}\label{sec:The structure of the paper}
The paper is organized as follows.
The background Section \ref{sec:Background} recalls briefly the construction of the exterior differential system \( P \) described above, whose solutions correspond to Pfaffian embeddings.
Next is a description of the facts on \( (2,3,5) \)-manifolds that will be used in the analysis, as well as the homogeneous space of isotropic \( 2 \)-planes.

The main part of the paper analyzes the coframed version \( \BP \) of \( P \), to determine the conditions necessary for solutions to exist.
After a precise statement of the characterization of embeddable \( (2,3,5) \)-manifolds (Subsection \ref{sec:The main theorem}), we start in Subsection \ref{sec:The coframed embedding PDE} with the construction of \( \BP \) and the description of its canonical linear Pfaffian coframing.
After doing so, the first integrability condition of the exterior differential system is easily described (\ref{sec:The first integrability condition}).

Subsection \ref{sec:Normalization after the first integrability condition} deals with the normalization of solutions to the first integrability condition.
The calculation there incidentally serves as a demonstration of the compatibility of symbols mentioned above.
After normalization, the tableau forms are described explicitly in terms of the Cartan geometric forms (\ref{sec:The tableau forms}), and the additional symbol relations from restricting to the first integrability are given (\ref{sec:The new symbol relations}).
The following two subsections \ref{sec:The additional integrability conditions} and \ref{sec:The remaining torsion absorption} deal with the remaining integrability conditions and torsion of the linear Pfaffian system.

Having absorbed all torsion, one can ask whether the differential system is in involution;
it is seen in Subsection \ref{sec:The failure of involutivity} that the equation is not involutive, so the following Subsection \ref{sec:Prolongation} deals with its prolongation, a process that enlarges the exterior differential system by \( 18 \) dimensions.
After prolongation, a combination of integrability conditions and symmetry reductions can be made to reduce the tableau to a trivial one, which is shown in Subsection \ref{sec:Higher reductions}.
This means in particular, that there is a set of obvious integrability conditions (depending only on the curvature functions of \( M \) as a \( (2,3,5) \)-geometry), which obstruct embeddings where nonzero.
On the other hand, if all of the integrability conditions vanish, then the differential system has been reduced to a Frobenius system, in which case the maximal leaves correspond to solutions of the original problem.
From the integrability conditions, several explicit restrictions on the curvature of embeddable \( M \) are given in Subsection \ref{sec:The curvature obstructions to embedding}.

In the final Section \ref{sec:Examples}, two examples of embeddable \( (2,3,5) \)-manifolds are given.
The first is, unsurprisingly, the flat model \( (2,3,5) \)-manifold.
It should be noted though, that this is already a non-trivial statement, because the \( (2,3,5) \)-symbol is not quite compatible for the application of Tanaka theory, and so a minor correction is needed.
The second example is an embeddable non-flat \( (2,3,5) \)-manifold.
Specifically, the example is a member of the \( 1 \)-parameter family of multiply-transitive Cartan geometries whose Cartan quartic has two double roots, symmetry type \( \mathbf{D.6}_{a} \) in the notation of \cite{TheCartantheoreticClassificationMultiplytransitive2022}.

\subsection{Acknowledgements}\label{sec:Acknowledgements}
The author would like to acknowledge valuable conversations on the topic of this paper with several people, including Jun-Muk Hwang, Tohru Morimoto, Dennis The, Robert Bryant, Igor Zelenko, Boris Doubrov, Andreas \v{C}ap, Michael Eastwood, Sean Curry, Omid Makhmali, Colleen Robles, and Paul-Andi Nagy.
For the calculations, the free open-source SageMath system was indispensable.

\subsection{AI statement}
At the time that this paper was written, AI was not yet effective at mathematics, so the ideas and writing were done without it.
There is one exception: on a more recent proofreading pass, Claude Opus 5 proposed to relate the unique embeddable two-double-root system uncovered here (see the example in Section \ref{sec:A non-flat example}) to known classifications.
This turns out to be straightforward, and Claude is capable of autonomously taking the results here and matching to the standard classifications.
The writeup of that statement in this paper is my own calculation, different from Claude's, but mathematically equivalent.
(The difference is primarily a matter of taste.)

\addtocontents{toc}{\protect\setcounter{tocdepth}{2}}

\section{Background}\label{sec:Background}
\subsection{The 1-jet contact system and the embedding PDE}\label{sec:The 1-jet contact system and the embedding PDE}
Given two smooth manifolds \( M \) and \( N \), denote by \( J^1 = J^{1}(M, N) \) the first jet space of maps \( M \) to \( N \), comprising the \( 1 \)-jets of germs of smooth maps \( M \to N \).
The space \( J^{1} \) is a bundle \( \pi_{0} \colon J^{1} \to M \times N \), with projection sending a \( 1 \)-jet \( j^{1}_{x}(f) \) based at \( x \) to its source and target, \( (x, f(x)) \in M \times N \).
The bundle \( J^{1} \) is naturally identified with the bundle \( \Hom(TM, TN) \) over \( M \times N \), with the identification on the fiber over each \( (x,y) \in M \times N \) given by the \emph{jet linearization} map
\begin{equation}\label{eq: jet linearization map}
   \begin{tikzcd}[row sep={8mm,between origins}]
    H_{(x,y)} \colon J^{1}_{(x,y)} \ar[r] & \Hom(T_{x} M, T_{y} N) \\
    \phantom{H_{(x,y)} \colon} j^{1}_{x}(f) \ar[r, mapsto] & \d f_{x} .
  \end{tikzcd} 
\end{equation}
Of course, \( J^{1} \) is also a bundle over \( M \) and \( N \), with respective projection maps
\[ J^{1} \xrightarrow{\pi_{M}} M, \qquad J^{1} \xrightarrow{\pi_{N}} N . \]

We recall the canonical \emph{contact form}\footnote{Here \emph{contact} is used in the older meaning, from which derives the special case now commonly used---a maximally non-integrable \( 1 \)-form on an odd-dimensional manifold. The contact form on \( J^{1}(\R^{n}, \R) \) is the prototypical example, and satisfies either definition.} on \( J^{1} \), defined as the \( TN \)-valued \( 1 \)-form
\[ \Theta_{u} = \d \pi_{N} - H_{\pi_{0}(u)}(u) \circ \d \pi_{M} \quad \mbox{ for } \quad u \in J^{1} . \]
Given a smooth section \( s \) of the bundle \( M \times N \to M \), the \emph{\( 1 \)-jet lift} of \( s \) is the section of \( J^{1} \to M \) sending each point \( x \in M \) to the \( 1 \)-jet of \( f = \pi_{N} \circ s \) at \( x \).
An arbitrary section of \( J^{1} \to M \) is the \( 1 \)-jet lift of some \( s \) if and only if it pulls back the contact form to zero.
This characterization of jet lifts is important, because smooth maps \( f \colon M \to N \) are in bijection with smooth sections \( s \) of \( M \times N \to M \), and these are in bijection with the \( 1 \)-jet lifts of smooth sections to \( J^{1} \to M \).

Now suppose that \( M \) and \( N \) have distributions \( D \) and respectively \( D' \).
For each fiber \( J^{1}_{(x,y)} \), let \( P_{(x,y)} \) consist of the injective linear maps that send \( D_{x} \) into \( D'_{y} \).
The union \( P \) of these is a subbundle of \( J^{1} \), and by construction, the sections of \( P \to M \) that pull back \( \Theta \) to zero are the 1-jet lifts corresponding to Pfaffian immersions \( (M,D) \to (N,D') \).
We call such a section an \emph{integral section} of \( P \to M \).
If we work locally, as we do here, these correspond exactly to Pfaffian embeddings.

So, the question now is: do there exist integral sections of \( P \to M \)?

\subsection{Curved (2,3,5)-manifolds}\label{sec:Curved (2,3,5)-manifolds}
Given a manifold \( M \) and a distribution \( D_{-1} \subset TM \), recall \cite{TanakaDifferentialSystemsGraded1970} that \( D_{-1} \) defines a \emph{regular differential system} if the inductively defined \emph{(weak) derived systems}\footnote{The underline here denotes the sheaves of sections of a vector bundle.}
\[ \underline{D}_{p-1} := [\underline{D}_{p}, \underline{D}_{-1}] \]
are of constant rank for each \( p < 0 \).
In this way one obtains from \( D_{-1} \) a \emph{filtered manifold} structure on \( M \), the filtration
\[ D_{-1} \subset D_{-2} \subset \ldots \subset TM \]
such that arbitrary sections \( X \) of \( D_{p} \) and \( Y \) of \( D_{q} \) have Lie bracket \( [X, Y] \) in \( D_{p+q} \).

It is a classical fact (Cartan \cite{CartanSystemesPfaffCinq1910}), that for a \( 5 \)-dimensional manifold \( M \) equipped with a \( 2 \)-dimensional plane distribution \( D_{-1} \), the generic behavior is for the first derived system to have rank \( 3 \) and the second derived system to have full rank \( 5 \), so that one has the filtration
\begin{equation}\label{eq: 235 filtration}
   D_{-1} \subset D_{-2} \subset D_{-3} = TM .
\end{equation}
Such a pair \( (M, D_{-1}) \) is called, by virtue of these ranks, a \( (2,3,5) \)-manifold.

Conversely, for any \( (2,3,5) \)-manifold, the derived distribution \( D_{-2} \) determines \( D_{-1} \) uniquely.
This can be seen using the \emph{Levi-bracket} on \( D_{-2} \), which is defined at each \( x \in M \) by the rule
\[ \begin{tikzcd}[row sep={8mm,between origins}]
   \Levi_{x, D_{-2}} \colon \Wedge^{2} (D_{-2})_{x} \ar[r] & T_{x}M / (D_{-2})_{x} \phantom{\mod{(D_{-2})}_{x}} \\
   \phantom{\Levi_{x, D_{-2}} \colon} X \w Y \ar[r, mapsto] & {[\tilde{X}, \tilde{Y}]_{x} \mod{D_{-2}}} ,
  \end{tikzcd} \]
where \( X, Y \) are vectors in \( (D_{-2})_{x} \) and \( \tilde{X}, \tilde{Y} \) are arbitrary local extensions to vector fields in \( D_{-2} \).
By assumption, \( \Levi_{x, D_{-2}} \) must be surjective, so that its kernel is \( 1 \)-dimensional and spanned by a decomposable \( 2 \)-vector.
This \( 2 \)-vector defines at each \( x \) a \( 2 \)-plane, so one recovers a \( 2 \)-plane distribution \( D'_{-1} \) from \( D_{-2} \).
On the other hand, \( \Wedge^{2} D_{-1} \) is in the kernel of \( \Levi_{x, D_{-2}} \), so \( D'_{-1} \) must equal \( D_{-1} \).
In fact, it is not difficult to check that a generic \( 3 \)-distribution on a \( 5 \)-manifold determines a \( (2,3,5) \)-structure in this manner on a dense open set.

The filtration \eqref{eq: 235 filtration} on \( M \) has associated graded bundle
\[ Gr = Gr_{-3} \oplus Gr_{-2} \oplus Gr_{-1} := D_{-3} / D_{-2} \oplus D_{-2} / D_{-1} \oplus D_{-1} , \]
and the Lie bracket of vector fields induces on each fiber of this bundle the structure of a nilpotent graded Lie algebra, using the same extension procedure as in the definition of the Levi-bracket (see \cite{TanakaDifferentialSystemsGraded1970} for details).
Thus we have, for example, a well defined map 
\[ Gr_{-1} \otimes Gr_{-2} \to Gr_{-3} , \]
which is necessarily surjective (else \( D_{-1} \) would not be bracket generating).
Since \( Gr_{-2} \) is \( 1 \)-dimensional, it follows that \( D_{-1} \) and \( TM/D_{-2} \) are canonically isomorphic, up to a scaling factor.

The Lie algebra structure on each fiber of \( Gr \) is isomorphic to a canonical Lie algebra \( \frg_{M,-} \).
Fix a two dimensional vector space \( U \) with basis \( e_{1}, e_{2} \), and let
\[ \frg_{M,-} = \frg_{-3} \oplus \frg_{-2} \oplus \frg_{-1} := U \oplus \R \oplus U \]
as vector spaces.
Define a basis 
\begin{equation}\label{eq: nilpotent 235 model basis}
   e_{1}, e_{2}, \quad e_{0}, \quad e_{1'}, e_{2'} 
\end{equation}
on \( \frg_{M,-} \), letting \( e_{0} = 1 \) and duplicating the basis on the second copy of \( U \) by marking with apostrophes.
For the Lie algebra structure on \( \frg_{M,-} \), define
\[ [e_{1'}, e_{2'}] = 2 e_{0}, \qquad
[e_{0}, e_{1'}] = 3 e_{1}, \qquad
[e_{0}, e_{2'}] = 3 e_{2}, \]
all other brackets vanishing.

In accordance with the basis \eqref{eq: nilpotent 235 model basis}, we will consistently use the indices \( a, b, c, \ldots \) ranging over the index set \( \{1, 2\} \), allowing also primed indices \( a', b', \ldots \) over the same range.
We will furthermore use the indices \( i, j, \ldots \) ranging over the index set \( \{0, 1', 2'\} \).

The prototypical example of a \( (2,3,5) \)-manifold is the maximally homogeneous model, \( G_{2} / P_{1} \), where \( G_{2} \) is the exceptional 14 dimensional Lie group and \( P_{1} \) a certain parabolic subgroup.
To describe \( P_{1} \), fix the standard representation embedding \( G_{2} \) into \( \GL(\R^{7}) \).
This embedding allows a compact description of the Lie algebra \( \frg_{2} \) via its Maurer-Cartan form,
\[ \left( \begin{matrix}
  -\zeta_{2} - 2 \zeta_{1} & \gamma^{0}_{2} & -\gamma^{0}_{1} & -\sqrt{2} \gamma & \gamma_{2} & \gamma_{1} & 0 \\
   \omega^{2'} & -\zeta_{2} - \zeta_{1} & -\gamma^{2}_{1} & -\sqrt{2} \gamma^{0}_{1} & \gamma & 0 & -\gamma_{1} \\
  -\omega^{1'} & -\gamma^{1}_{2} &-\zeta_{1} & -\sqrt{2} \gamma^{0}_{2} & 0 & -\gamma & -\gamma_{2} \\
  \sqrt{2} \omega^{0} &-\sqrt{2} \omega^{1'} & -\sqrt{2} \omega^{2'} & 0 & \sqrt{2} \gamma^{0}_{2} & \sqrt{2} \gamma^{0}_{1} & \sqrt{2} \gamma \\
  \theta^{2} & -\omega^{0} & 0 & \sqrt{2} \omega^{2'} & \zeta_{1} & \gamma^{2}_{1} &\gamma^{0}_{1} \\
  \theta^{1} & 0 & \omega^{0} & \sqrt{2} \omega^{1'} & \gamma^{1}_{2} &\zeta_{2} + \zeta_{1} & -\gamma^{0}_{2} \\
  0 & -\theta^{1} & -\theta^{2} &-\sqrt{2} \omega^{0} & \omega^{1'} & -\omega^{2'} & \zeta_{2} + 2 \zeta_{1}
\end{matrix} \right) , \]
with the component \( 1 \)-forms 
\[ \theta^{1}, \theta^{2}, \omega^{0}, \omega^{1'}, \omega^{2'}, \gamma^{1}_{2}, \zeta_{1}, \zeta_{2}, \gamma^{2}_{1}, \gamma^{0}_{1}, \gamma^{0}_{2}, \gamma, \gamma_{1}, \gamma_{2} \]
defining a parallelism of \( G_{2} \).
Denote the dual left invariant vector fields to these by
\begin{equation}\label{eq: G2 vf basis}
   X_{\theta^{1}}, X_{\theta^{2}}, X_{\omega^{0}}, X_{\omega^{1'}}, X_{\omega^{2'}}, X_{\gamma^{1}_{2}}, X_{\zeta_{1}}, X_{\zeta_{2}}, X_{\gamma^{2}_{1}}, X_{\gamma^{0}_{1}}, X_{\gamma^{0}_{2}}, X_{\gamma}, X_{\gamma_{1}}, X_{\gamma_{2}} . 
\end{equation}
In this embedding, \( P_{1} \) is the subgroup of \( G_{2} \) that preserves the line through the lowest weight vector \( (1,0,0,0,0,0,0)^{t} \).

At this point we turn to graded structures, so it is necessary to change notation, denoting \( G_{M} := G_{2} \) and \( G^{0}_{M} := P_{1} \) as well as their Lie algebras \( \frg_{M} := \frg_{2} \) and \( \frg^{0}_{M} := \mathfrak{p}_{1} \).
By standard theory, e.g. \cite{Cap.SlovakParabolicGeometries2009} Chapter 3, a choice of parabolic subalgebra is the same as a choice of grading on \( \frg_{M} \) compatible with the Lie bracket.
In this case the grading is
  \[ \frg_{M} = \frg_{M,-3} \oplus \frg_{M,-2} \oplus \frg_{M,-1} \oplus \frg_{M,0} \oplus \frg_{M,1} \oplus \frg_{M,2} \oplus \frg_{M,3} \]
where the sum of degree negative gradeds
\[ \frg_{M,-3} \oplus \frg_{M,-2} \oplus \frg_{M,-1} = \{X_{\theta^{1}}, X_{\theta^{2}} \} \oplus \{ X_{\omega^{0}} \} \oplus \{ X_{\omega^{1'}}, X_{\omega^{2'}} \} \]
is visibly isomorphic to \( \frg_{M,-} \) as a Lie algebra,
the degree zero graded
\[ \frg_{M,0} = \{ X_{{\gamma}^{1}_{2}}, X_{{\zeta}_{1}}, X_{{\zeta}_{2}}, X_{{\gamma}^{2}_{1}} \} \]
is isomorphic to \( \gl(U) \),
and the sum of degree positive gradeds
\[ \frg_{M,+} = \frg_{M,1} \oplus \frg_{M,2} \oplus \frg_{M,3} = \{ X_{{\gamma}^{0}_{1}}, X_{{\gamma}^{0}_{2}} \} \oplus \{ X_{{\gamma}} \} \oplus \{ X_{{\gamma}_{1}}, X_{{\gamma}_{2}} \} \]
is isomorphic to the dual of \( \frg_{M,-} \).

The grading on \( \frg_{M} \) induces a filtration, the \emph{upper degree} filtration, by
\[ \frg^{p}_{M} = \bigoplus_{q \ge p} \frg_{M,q} \quad \mbox{ for } \quad p \in \Z . \]
Note in particular that there is a canonical isomorphism
\[ \frg_{M,-} \cong \frg_{M} / \frg^{0}_{M} . \]
The notation is consistent with the parabolic subalgebra \( \frg^{0}_{M} \), which is simply described as the non-negative graded component of \( \frg_{M} \).
Likewise, \( G^{0}_{M} \) can be characterized as the subgroup of elements in \( G_{M} \) whose adjoint action preserves the upper degree filtration of \( \frg_{M} \).
More generally, for \( p \ge 0 \), are the subgroups \( G^{p}_{M} \) of elements \( g \in G^{0}_{M} \) such that \( \Ad_{g}(\frg^{q}_{M}) \subseteq \frg^{q+p}_{M} \) for all \( q \in \Z \).
It will be useful to note that the exponential map of \( \frg_{M} \) is a diffeomorphism when restricted to \( \frg^{1}_{M} \), because \( \frg^{1}_{M} \) is nilpotent (see for example \cite{Cap.SlovakParabolicGeometries2009} Theorem 3.1.3).

As just noted, the flat \( (2,3,5) \)-manifold gives a good model for general \( (2,3,5) \)-manifolds.
Precisely, \( (2,3,5) \)-manifolds are in bijection with \emph{regular}, \emph{normal} \emph{Cartan geometries} \( (\mathcal{B}_{M}, \Omega_{M}) \) modeled on \( G_{M} / G^{0}_{M} \).
By definition, \( \mathcal{B}_{M} \) is a \( G^{0}_{M} \)-principal bundle over \( M \) and \( \Omega_{M} \) a \( \frg_{M} \)-valued \emph{Cartan connection} form on \( \mathcal{B}_{M} \).
Recall that for \( \Omega_{M} \) to be a Cartan connection means
\begin{itemize}
\item \( \Omega_{M} \) gives an isomorphism of each tangent space of \( \B_{M} \) to \( \frg_{M} \),
\item for each point \( p \in \B_{M} \) and vector \( X \in \frg^{0}_{M} \), the isomorphism \( \Omega_{M} \) reproduces the vector \( \tilde{X} \in T_{p} \B_{M} \) induced by \( X \),
\[ \Omega_{M}(\tilde{X}) = X , \]
\item and \( \Omega_{M} \) is \( G^{0}_{M} \)-equivariant,
\[ R^{*}_{g}\Omega_{M} = \Ad_{g^{-1}} \circ \Omega_{M} \quad \mbox{ for } \quad g \in G^{0}_{M} . \]
\end{itemize}
Refer to \cite{Cap.SlovakParabolicGeometries2009} for more details on this correspondence and Cartan geometries in general, and to \cite{TheCartantheoreticClassificationMultiplytransitive2022} for the specific case of \( (2,3,5) \) and some history.

Using the basis of \( \mathfrak{g}_{2} \) fixed in display \eqref{eq: G2 vf basis} , the parallelism \( \Omega_{M} \) gives rise to a canonical tautological coframing on \( \mathcal{B}_{M} \), denoted again by
\[ \theta^{1}, \theta^{2}, \omega^{0}, \omega^{1'}, \omega^{2'}, \quad \gamma^{1}_{2}, \zeta_{1}, \zeta_{2}, \gamma^{2}_{1}, \gamma^{0}_{1}, \gamma^{0}_{2}, \gamma, \gamma_{1}, \gamma_{2} . \]
In case \( M \) is flat, this coframing is of course the Maurer-Cartan coframing described above.
We note for the sake of comparison that the corresponding forms in \cite{TheCartantheoreticClassificationMultiplytransitive2022} are
\[ \theta_{32}, \theta_{31}, \theta_{21}, -\theta_{11}, \theta_{10}, \quad -\theta_{01}, -\zeta_{1}, -\zeta_{2}, -\pi_{01}, -\pi_{11}, \pi_{10}, -\pi_{21}, \pi_{32}, \pi_{31} . \]

Recall that any Cartan geometry has an associated first order \( G \)-structure \( \coframeBundle_{M} \), as in \cite{Cap.SlovakParabolicGeometries2009}.
In the case at hand, \( \coframeBundle_{M} \) is the \( G^{0}_{M} / G^{3}_{M} \)-structure obtained as the quotient of \( \B_{M} \) by the action of \( G^{3}_{M} \),
due to the fact that \( G^{3}_{M} \) is the kernel of the adjoint action modulo \( \frg^{0}_{M} \),
\[ \begin{tikzcd}[row sep={8mm,between origins}]
    G^{0}_{M} \ar[r] & \GL(\frg_{M} / \frg^{0}_{M}) \\
    A \ar[r, mapsto] & \left(X \mapsto Ad_{A} X \mod{\frg^{0}_{M}}\right)
  \end{tikzcd} . \]
Note in particular the embedding of \( G^{0}_{M} / G^{3}_{M} \) as a subgroup of \( \GL(\frg_{M,-}) \cong \GL(\frg_{M} / \frg^{0}_{M}) \), and the corresponding embedding \( \frg_{M,0} \oplus \frg_{M,1} \oplus \frg_{M,2} \to \gl(\frg_{M,-}) \).
The soldering form of \( \coframeBundle_{M} \) is obtained from the \( \frg_{M,-} \)-component \( \Omega_{M,-} \) of \( \Omega_{M} \), which satisfies the equation
\[ \Omega_{M,-}(X) = \Omega_{M}(X) \mod{\frg^{0}_{M}} , \]
again using the identification \( \frg_{M,-} \cong \frg_{M} / \frg^{0}_{M} \).
The right hand side is manifestly \( G^{3}_{M} \)-invariant, so this form descends to \( \coframeBundle_{M} \).
Otherwise put, for each point \( p \in \B_{M} \), the Cartan connection modulo \( \frg^{0}_{M} \) defines a linear frame
\[ \begin{tikzcd}
    \varphi_{M,p} \colon T_{x} M \ar[r, "\cong"] & \frg_{M,-} .
  \end{tikzcd} \]
It is not difficult to check that \( \varphi_{M,p\cdot g} = \Ad_{g^{-1}} \circ \varphi_{M,p} \) for \( g \in G^{0}_{M} \).

Due to the construction, we can safely regard \( \Omega_{M,-} \) as a tautological coframing form for \( M \).
For this reason, it will be convenient to use the basis \( e_{a}, e_{i} \) to decompose \( \Omega_{M,-} \) into \( 1 \)-forms \( \Omega^{a}, \Omega^{i} \).

The Cartan form also defines a connection form on \( \coframeBundle_{M} \).
Indeed, the \( \frg_{M,0} \oplus \frg_{M,1} \oplus \frg_{M,2} \)-component of \( \Omega_{M} \) is \( G^{3}_{M} \) invariant, so descends to a form on \( \coframeBundle_{M} \).
This quotient form is \( G^{0}_{M} / G^{3}_{M} \)-equivariant and complementary to \( \Omega_{M,-} \), so defines a connection form.
Because \( \frg_{M,0} \oplus \frg_{M,1} \oplus \frg_{M,2} \) embeds into \( \gl(\frg_{M,-}) \), we may regard the connection form as a \( \gl(\frg_{M,-}) \)-valued form.
In terms of the basis \( e_{a}, e_{i} \), it has components \( \Omega^{a}_{b}, \Omega^{i}_{a}, \Omega^{i}_{j} \), with the \( \Omega^{a}_{i} \) component vanishing identically.

The first order structure equations of \( \B_{M} \) can be found in \cite{TheCartantheoreticClassificationMultiplytransitive2022} Appendix A, and reflect the \( G^{0}_{M} / G^{3}_{M} \)-structure just described,
\begin{equation}\label{eq: 235 first order structure}
  \d \left(
    \begin{matrix}
      \theta^{1} \\ \theta^{2} \\ \omega^{0} \\ \omega^{1'} \\ \omega^{2'}
    \end{matrix} \right)
  ={}
  -\left(
    \begin{matrix}
      3\zeta_{1} + 2\zeta_{2} & \gamma^{1}_{2} & 0 & 0 & 0 \\
      \gamma^{2}_{1} & 3\zeta_{1} + \zeta_{2} & 0 & 0 & 0 \\
      \gamma^{0}_{1} & \gamma^{0}_{2} & 2\zeta_{1} + \zeta_{2} & 0 & 0 \\
      \gamma & 0 & 2\gamma^{0}_{2} & \zeta_{1} + \zeta_{2} & \gamma^{1}_{2} \\
      0 & \gamma & -2 \gamma^{0}_{1} & \gamma^{2}_{1} & \zeta_{1}
    \end{matrix} \right)\!\! \w \!\!
  \left( \begin{matrix}
      \theta^{1} \\ \theta^{2} \\ \omega^{0} \\ \omega^{1'} \\ \omega^{2'}
    \end{matrix} \right)
  + \left(
    \begin{matrix}
      3 \omega^{0} \w \omega^{1'} \\ 3 \omega^{0} \w \omega^{2'}
      \\ 2\omega^{1'} \w \omega^{2'} \\ 0 \\ 0
    \end{matrix} \right) .
\end{equation}
We have in particular that the connection form is
\[
\left(
 \begin{matrix}
   \Omega^{a}_{b} & 0 & 0 \\
   \Omega^{0}_{b} & \Omega^{0}_{0} & 0 \\
   \Omega^{a'}_{b} & \Omega^{a'}_{0} & \Omega^{a'}_{b'}
 \end{matrix} \right)
= \left(
 \begin{matrix}
   3\zeta_{1} + 2\zeta_{2} & \gamma^{1}_{2} & 0 & 0 & 0 \\
   \gamma^{2}_{1} & 3\zeta_{1} + \zeta_{2} & 0 & 0 & 0 \\
   \gamma^{0}_{1} & \gamma^{0}_{2} & 2\zeta_{1} + \zeta_{2} & 0 & 0 \\
   \gamma & 0 & 2\gamma^{0}_{2} & \zeta_{1} + \zeta_{2} & \gamma^{1}_{2} \\
   0 & \gamma & -2 \gamma^{0}_{1} & \gamma^{2}_{1} & \zeta_{1}
 \end{matrix} \right) .
\]
The \emph{torsion functions} of these structure equations are the \( \frg_{M,-} \otimes \Wedge^{2} \frg^{\vee}_{M,-} \)-valued functions \( T \) such that
\[ \d\Omega^{\alpha} = - \Omega^{\alpha}_{\beta}\w\Omega^{\beta} + T^{\alpha}_{\beta,\gamma} \Omega^{\beta} \w \Omega^{\gamma} , \]
temporarily letting \( \alpha, \beta, \gamma = 1, 2, 0, 1', 2' \),
so that
\begin{equation}\label{eq: M torsion functions}
   T^{\thindex{1}}_{0,1'} = 3, \quad T^{\thindex{2}}_{0,2'} = 3, \quad T^{0}_{1',2'} = 2 ,  
\end{equation}
with all other components zero.
Here we have employed the Einstein summation convention, and will continue to do so where appropriate.

The higher order structure equations, also found in \cite{TheCartantheoreticClassificationMultiplytransitive2022}, are
\begin{equation}\label{eq: 235 higher order structure}
  \begin{aligned}
    \d \gamma^{1}_{2} & = -\zeta_{2}\w\gamma^{1}_{2} {} + \gamma_{2}\w\theta^{1} + 3\gamma^{0}_{2}\w\omega^{1'}  + K_{\gamma^{1}_{2}}
    \\ \d \zeta_{1}  & = \gamma^{1}_{2}\w\gamma^{2}_{1} {} + \gamma_{2}\w\theta^{2} - \gamma\w\omega^{0} - \gamma^{0}_{1}\w\omega^{1'} + 2\gamma^{0}_{2}\w\omega^{2'}  + K_{\zeta_{1}}
    \\ \d \zeta_{2}  & = -2\gamma^{1}_{2}\w\gamma^{2}_{1} {} + \gamma_{1}\w\theta^{1} - \gamma_{2}\w\theta^{2} + 3\gamma^{0}_{1}\w\omega^{1'} - 3\gamma^{0}_{2}\w\omega^{2'}  + K_{\zeta_{2}}
    \\ \d \gamma^{2}_{1} & =  \zeta_{2}\w\gamma^{2}_{1} {} + \gamma_{1}\w\theta^{2} + 3\gamma^{0}_{1}\w\omega^{2'}  + K_{\gamma^{2}_{1}}
    \\ \d \gamma^{0}_{1} & =  \gamma^{2}_{1}\w\gamma^{0}_{2} + \zeta_{2}\w\gamma^{0}_{1} + \zeta_{1}\w\gamma^{0}_{1} {} + \gamma_{1}\w\omega^{0} + 2\gamma\w\omega^{2'}  + K_{\gamma^{0}_{1}}
    \\ \d \gamma^{0}_{2} & = \gamma^{1}_{2}\w\gamma^{0}_{1} + \zeta_{1}\w\gamma^{0}_{2} {} + \gamma_{2}\w\omega^{0} - 2\gamma\w\omega^{1'}  + K_{\gamma^{0}_{2}}
    \\ \d \gamma & =  -2\gamma^{0}_{2}\w\gamma^{0}_{1} + \zeta_{2}\w\gamma + 2\zeta_{1}\w\gamma {} + \gamma_{1}\w\omega^{1'} + \gamma_{2}\w\omega^{2'}  + K_{\gamma}
    \\ \d \gamma_{1} & =  -3\gamma^{0}_{1}\w\gamma + \gamma^{2}_{1}\w\gamma_{2} + 2\zeta_{2}\w\gamma_{1} + 3\zeta_{1}\w\gamma_{1} {}  + K_{\gamma_{1}}
    \\ \d \gamma_{2} & =  -3\gamma^{0}_{2}\w\gamma + \zeta_{2}\w\gamma_{2} + 3\zeta_{1}\w\gamma_{2} + \gamma^{1}_{2}\w\gamma_{1} {}  + K_{\gamma_{2}}
  \end{aligned}
\end{equation}
where the \( 2 \)-forms \( K_{\gamma^{1}_{2}}, \ldots, K_{\gamma_{2}} \) comprise the curvature of the Cartan geometry, explicitly given as
\begin{align*}
  K_{\gamma^{1}_{2}} ={}& A_{1} \theta^{2}\w\omega^{2'} + A_{2}(\theta^{1}\w\omega^{2'} + \theta^{2}\w\omega^{1'}) + A_{3} \theta^{1}\w\omega^{1'} - 2 B_{1} \theta^{2}\w\omega^{0} - 2 B_{2} \theta^{1}\w\omega^{0} + C_{1} \theta^{1}\w\theta^{2}
  \\ K_{\zeta_{1}} ={}& - A_{2} \theta^{2}\w\omega^{2'} - A_{3} (\theta^{1}\w\omega^{2'} + \theta^{2}\w\omega^{1'}) - A_{4} \theta^{1}\w\omega^{1'} + 2 B_{2} \theta^{2}\w\omega^{0} + 2 B_{3} \theta^{1}\w\omega^{0} -C_{2} \theta^{1}\w\theta^{2}
  \\ K_{\zeta_{2}} ={}& 2 A_{2} \theta^{2}\w\omega^{2'} + 2 A_{3} (\theta^{1}\w\omega^{2'} + \theta^{2}\w\omega^{1'}) + 2 A_{4} \theta^{1}\w\omega^{1'} - 4 B_{2} \theta^{2}\w\omega^{0} - 4 B_{3} \theta^{1}\w\omega^{0} + 2 C_{2} \theta^{1}\w\theta^{2}
  \\ K_{\gamma^{2}_{1}} ={}& - A_{3} \theta^{2}\w\omega^{2'} - A_{4} (\theta^{1}\w\omega^{2'} + \theta^{2}\w\omega^{1'}) - A_{5} \theta^{1}\w\omega^{1'} + 2 B_{3} \theta^{2}\w\omega^{0} + 2 B_{4} \theta^{1}\w\omega^{0} -C_{3} \theta^{1}\w\theta^{2}
  \\ K_{\gamma^{0}_{1}} ={}& - B_{2} \theta^{2}\w\omega^{2'} - B_{3} (\theta^{1}\w\omega^{2'} + \theta^{2}\w\omega^{1'}) - B_{4} \theta^{1}\w\omega^{1'} + 2 C_{2} \theta^{2}\w\omega^{0} + 2 C_{3} \theta^{1}\w\omega^{0} -D_{2} \theta^{1}\w\theta^{2}
  \\ K_{\gamma^{0}_{2}} ={}& - B_{1} \theta^{2}\w\omega^{2'} - B_{2} (\theta^{1}\w\omega^{2'} + \theta^{2}\w\omega^{1'}) - B_{3} \theta^{1}\w\omega^{1'} + 2 C_{1} \theta^{2}\w\omega^{0} + 2 C_{2} \theta^{1}\w\omega^{0} -D_{1} \theta^{1}\w\theta^{2}
  \\ K_{\gamma} ={}& - C_{1} \theta^{2}\w\omega^{2'} - C_{2} (\theta^{1}\w\omega^{2'} + \theta^{2}\w\omega^{1'}) - C_{3} \theta^{1}\w\omega^{1'} + 2 D_{1} \theta^{2}\w\omega^{0} + 2 D_{2} \theta^{1}\w\omega^{0} -E \theta^{1}\w\theta^{2}
  \\ K_{\gamma_{1}} ={}& \tfrac{2}{3} D_{1} \theta^{2}\w\omega^{2'} + \tfrac{1}{3} D_{2} (\theta^{1}\w\omega^{2'} + \theta^{2}\w\omega^{1'}) - E \theta^{2}\w\omega^{0} - \tilde{D}_{2} \theta^{2}\w\omega^{2'} - \tilde{D}_{3} (\theta^{1}\w\omega^{2'} + \theta^{2}\w\omega^{1'}) \\  & - \tilde{D}_{4} \theta^{1}\w\omega^{1'} + 2 \tilde{E}_{2} \theta^{2}\w\omega^{0} + 2 \tilde{E}_{3} \theta^{1}\w\omega^{0} -\tilde{F}_{2} \theta^{1}\w\theta^{2}
  \\ K_{\gamma_{2}} ={}& - \tfrac{1}{3} D_{1} (\theta^{1}\w\omega^{2'} + \theta^{2}\w\omega^{1'}) - \tfrac{2}{3} D_{2} \theta^{1}\w\omega^{1'} + E \theta^{1}\w\omega^{0} - \tilde{D}_{1} \theta^{2}\w\omega^{2'} - \tilde{D}_{2} (\theta^{1}\w\omega^{2'} + \theta^{2}\w\omega^{1'}) \\
  & - \tilde{D}_{3} \theta^{1}\w\omega^{1'} + 2 \tilde{E}_{1} \theta^{2}\w\omega^{0} + 2 \tilde{E}_{2} \theta^{1}\w\omega^{0} -\tilde{F}_{1} \theta^{1}\w\theta^{2}
  \end{align*}
for \emph{curvature functions}
\begin{equation*}
  \begin{gathered}
    A_{1}, A_{2}, A_{3}, A_{4}, A_{5}, \enspace B_{1}, B_{2}, B_{3}, B_{4}, \enspace C_{1}, C_{2}, C_{3}, \enspace D_{1}, D_{2}, \enspace E,
    \\ \tilde{D}_{1}, \tilde{D}_{2}, \tilde{D}_{3}, \tilde{D}_{4}, \enspace \tilde{E}_{1}, \tilde{E}_{2}, \tilde{E}_{3}, \enspace \tilde{F}_{1}, \tilde{F}_{2} .
  \end{gathered}
\end{equation*}
These curvature functions are equivariant functions on \( \mathcal{B}_{M} \) taking values in various \( G^{0}_{M} \)-modules.
Accordingly, the differential of each curvature function is the sum of an equivariant component and a semibasic (linear combination of the \( \theta^{a}, \omega^{i} \)) component.
For example, one can check that
\[ \d A_{1} = 4 A_{2} \gamma^{1}_{2} + 4 A_{1} \zeta_{1} + A_{1;1} \theta^{1} + A_{1;2} \theta^{2} + A_{1;0} \omega^{0} + A_{1;1'} \omega^{1'} + A_{1;2'} \omega^{2'} \]
where the functions \( A_{1;a}, A_{1;i} \) on \( \mathcal{B}_{M} \) are the \emph{semibasic derivatives} of \( A_{1} \).
We will denote all other semibasic derivatives of the curvature functions likewise, referring to \cite{TheCartantheoreticClassificationMultiplytransitive2022} Appendix A for the full table.
The same discussion holds true for the semibasic derivatives of semibasic derivatives, and so on.
We will use similar notation for these, so that for example
\[ \d A_{5;0} = \mbox{ (equivariant terms) } + A_{5;0;1} \theta^{1} + A_{5;0;2} \theta^{2} + A_{5;0;0} \omega^{0} + A_{5;0;1'} \omega^{1'} + A_{5;0;2'} \omega^{2'} . \]
In fact, \( A_{5;0;1'} \) is the only semibasic second derivative that we will need to refer to.

As a final note, there are many functional relations between these higher derivatives, as can be seen from the requirement that \( \d^{2}\Gamma = 0 \) for each \( \Gamma \) in the coframing.
For example, the \( \theta^{2}\w\omega^{0}\w\omega^{1'} \) coefficient of \( \d^{2} \gamma^{2}_{1} \) is calculated to be 
\begin{equation}\label{eq: higher curvature relation example}
   -3 C_{3} + A_{4;0}  +2B_{3;1'} ,
\end{equation}
so determines that as a relation.
Such relations will be needed in Subsection \ref{sec:The curvature obstructions to embedding}, but there are more relations than is worthwhile to list here, especially as they can easily be determined as needed.
\subsection{Flat (4,7)-manifolds}\label{sec:Flat (4,7)-manifolds}
We now turn to the space of isotropic \( 2 \)-planes in a \( 6 \)-dimensional symplectic vector space.
To this end, consider \( \R^{6} \) as a vector space equipped with the standard symplectic structure \( \alpha \), and the space
\( N := \operatorname{Iso}_{2}(\R^{6}) \)
of isotropic \( 2 \)-planes in \( \R^{6} \).
The group \( G_{N} := \Sp(\R^{6}) \) acts transitively on \( N \), and the subgroup fixing a given point of \( N \) is a parabolic subgroup that we denote \( G^{0}_{N} \).
As such, \( N \) is a homogeneous space, isomorphic to \( G_{N} / G^{0}_{N} \).
The Maurer-Cartan form \( \Omega_{N} \) of \( G_{N} \) determines a flat Cartan connection on the \( G^{0}_{N} \)-principal bundle \( G_{N} \to N \).

To describe the Maurer-Cartan form, fix a symplectic basis of \( \R^{6} \),
\begin{equation}\label{eq: symplectic basis}
 e^{1}, e^{2}, e_{3}, e_{1}, e_{2}, e_{3'}
\end{equation}
such that
\[ \alpha(e^{1}, e_{1}) = \alpha(e^{2}, e_{2}) = \alpha(e_{3}, e_{3'}) = 1 , \]
all other pairs evaluating to zero.
Relative this basis is an embedding of \( \Sp(\R^{6}) \) into \( \GL(\R^{6}) \), in terms of which the Maurer-Cartan form is realized as
\[
\left(
  \begin{matrix}
    -\eta^{1}_{1} & - \eta^{2}_{1} & - \eta_{\tendex{1}{3}} & \eta_{11} & \eta_{12} & \eta_{\tendex{1}{3'}} \\
    -\eta^{1}_{2} & - \eta^{2}_{2} & - \eta_{\tendex{2}{3}} & \eta_{12} & \eta_{22} & \eta_{\tendex{2}{3'}} \\
    \varpi^{\tendex{1}{3}} & \varpi^{\tendex{2}{3}} & \eta^{3}_{3} & \eta_{\tendex{1}{3'}} & \eta_{\tendex{2}{3'}} & \eta^{3}_{3'} \\
    \vartheta^{\symdex{1}{1}} & \vartheta^{\symdex{1}{2}}  & \varpi^{\tendex{1}{3'}} & \eta^{1}_{1} & \eta^{1}_{2} & - \varpi^{\tendex{1}{3}} \\
    \vartheta^{\symdex{1}{2}} & \vartheta^{\symdex{2}{2}}  & \varpi^{\tendex{2}{3'}} & \eta^{2}_{1} & \eta^{2}_{2} & - \varpi^{\tendex{2}{3}} \\
    \varpi^{\tendex{1}{3'}} & \varpi^{\tendex{2}{3'}} & \eta^{3'}_{3} & \eta_{\tendex{1}{3}} & \eta_{\tendex{2}{3}} & - \eta^{3}_{3}
  \end{matrix} \right) ,
\]
where the component \( 1 \)-forms
\begin{align*}
\vartheta^{11}, \vartheta^{12}, \vartheta^{22}, \varpi^{\tendex{1}{3}}, \varpi^{\tendex{1}{3'}}, \varpi^{\tendex{2}{3}}, \varpi^{\tendex{2}{3'}}, \qquad & \eta^{1}_{1}, \eta^{1}_{2}, \eta^{2}_{1}, \eta^{2}_{2}, \eta^{3}_{3}, \eta^{3}_{3'}, \eta^{3'}_{3}, \\
&\eta_{\tendex{1}{3}}, \eta_{\tendex{1}{3'}}, \eta_{\tendex{2}{3}}, \eta_{\tendex{2}{3'}}, \eta_{11}, \eta_{12}, \eta_{22}
\end{align*}
define a complete coframing of \( \mathcal{B}_{N} \).
Let the corresponding dual basis be given by
\begin{align*}
e_{\symdex{1}{1}}, e_{\symdex{1}{2}}, e_{\symdex{2}{2}}, e_{\tendex{1}{3}}, e_{\tendex{1}{3'}}, e_{\tendex{2}{3}}, e_{\tendex{2}{3'}},
\quad & X_{\eta^{1}_{1}}, X_{\eta^{1}_{2}}, X_{\eta^{2}_{1}}, X_{\eta^{2}_{2}}, X_{\eta^{3}_{3}}, X_{\eta^{3}_{3'}}, X_{\eta^{3'}_{3}}, \\
& X_{\eta_{\tendex{1}{3}}}, X_{\eta_{\tendex{1}{3'}}}, X_{\eta_{\tendex{2}{3}}}, X_{\eta_{\tendex{2}{3'}}}, X_{\eta_{11}}, X_{\eta_{12}}, X_{\eta_{22}} .
\end{align*}

The grading of the Lie algebra \( \frg_{N} = \sp(6) \) compatible with the parabolic subgroup \( G^{0}_{N} \) is
\[ \frg_{N} = \frg_{N,-2} \oplus \frg_{N,-1} \oplus \frg_{N,0} \oplus \frg_{N,1} \oplus \frg_{N,2} , \]
where
\begin{align*}
 \frg_{N,-} = \frg_{N,-2} \oplus \frg_{N,-1}
 & = \{ e_{\symdex{1}{1}}, e_{\symdex{1}{2}}, e_{\symdex{2}{2}} \}
 \oplus \{ e_{\tendex{1}{3}}, e_{\tendex{1}{3'}}, e_{\tendex{2}{3}}, e_{\tendex{2}{3'}} \} ,
 \\ \frg_{N,0} & = \{ X_{\eta^{1}_{1}}, X_{\eta^{1}_{2}}, X_{\eta^{2}_{1}}, X_{\eta^{2}_{2}}, X_{\eta^{3}_{3}}, X_{\eta^{3}_{3'}}, X_{\eta^{3'}_{3}} \} , 
 \\ \frg_{N,1} \oplus \frg_{N,2}
 & = \{ X_{\eta_{\tendex{1}{3}}}, X_{\eta_{\tendex{1}{3'}}}, X_{\eta_{\tendex{2}{3}}}, X_{\eta_{\tendex{2}{3'}}} \}
\oplus \{ X_{\eta_{11}}, X_{\eta_{12}}, X_{\eta_{22}} \} . 
\end{align*}
The Lie algebra structure of \( \frg_{N,-} \) can be read directly from the Maurer-Cartan form, and is given by
\begin{align*}
[e_{\tendex{1}{3}}, e_{\tendex{1}{3'}}] = 2 e_{\symdex{1}{1}} & \qquad [e_{\tendex{1}{3}}, e_{\tendex{2}{3'}}] =  e_{\symdex{1}{2}} \\
[e_{\tendex{2}{3}}, e_{\tendex{2}{3'}}] = 2 e_{\symdex{2}{2}} & \qquad [e_{\tendex{2}{3}}, e_{\tendex{1}{3'}}] =  e_{\symdex{1}{2}} ,
\end{align*}
all other brackets vanishing.

We recall the discussion in Subsection \ref{sec:Curved (2,3,5)-manifolds} on the \( G \)-structure induced by a Cartan geometry.
In this case, the \( G \)-structure \( \coframeBundle_{N} \) has structure group \( G^{0}_{N} / G^{2}_{N} \).
The \( \frg_{N,-} \)-component of \( \Omega_{N} \) descends to the soldering form of \( \coframeBundle_{N} \), and the \( \frg_{N,0} \oplus \frg_{N,1} \)-component descends to a connection form on \( \coframeBundle_{N} \).
In particular, each point \( p \in \B_{N} \) over \( y \in N \) defines a linear frame
\[ \begin{tikzcd}
    \varphi_{N,p} \colon T_{y} N \ar[r, "\cong"] & \frg_{N,-} ,
  \end{tikzcd} \]
which satisfies \( \varphi_{N,p\cdot g} = \Ad_{g^{-1}} \circ \varphi_{N,p} \) for \( g \in G^{0}_{N} \).
Furthermore, with respect to the left multiplication action of \( G_{N} \) on \( \B_{N} = G_{N} \), is the equation
\( \varphi_{N,g\cdot p} = \varphi_{N,p} \circ \d L_{g^{-1}} \) for any \( g \in G_{N} \), which follows from the definition and the fact that the Cartan connection is simply the left invariant Maurer-Cartan form.

The structure equations follow simply from the Maurer-Cartan equations.
At first order,
  \begin{equation}\label{eq: N structure equations}
    \d \left( \begin{matrix}
    \vartheta^{\symdex{a}{b}} \\ \varpi^{\tendex{a}{3}} \\ \varpi^{\tendex{a}{3'}}
  \end{matrix} \right)
= -\left( \begin{matrix}
    \eta^{\symdex{a}{b}}_{\symdex{c}{d}} & 0 & 0 \\
    \eta^{\tendex{a}{3}}_{\symdex{c}{d}} & \eta^{\tendex{a}{3}}_{\tendex{c}{3}} & \eta^{\tendex{a}{3}}_{\tendex{c}{3'}} \\
    \eta^{\tendex{a}{3'}}_{\symdex{c}{d}} & \eta^{\tendex{a}{3'}}_{\tendex{c}{3}} & \eta^{\tendex{a}{3'}}_{\tendex{c}{3'}}
  \end{matrix} \right) \w
\left( \begin{matrix}
    \vartheta^{\symdex{c}{d}} \\ \varpi^{\tendex{c}{3}} \\ \varpi^{\tendex{c}{3'}}
  \end{matrix} \right)
+ \left( \begin{matrix}
    T^{\vartheta^{\symdex{a}{b}}} \\ T^{\tendex{a}{3}} \\ T^{\tendex{a}{3'}}
    \end{matrix} \right) ,
\end{equation}
where the component pieces of the connection form are 
\[ \left(
    \begin{matrix}
      \Omega^{ab}_{cd}
    \end{matrix} \right) =
  \left( \begin{matrix}
      \eta^{\symdex{1}{1}}_{\symdex{1}{1}} & \eta^{\symdex{1}{1}}_{\symdex{1}{2}} & \eta^{\symdex{1}{1}}_{\symdex{2}{2}} \\
      \eta^{\symdex{1}{2}}_{\symdex{1}{1}} & \eta^{\symdex{1}{2}}_{\symdex{1}{2}} & \eta^{\symdex{1}{2}}_{\symdex{2}{2}} \\
      \eta^{\symdex{2}{2}}_{\symdex{1}{1}} & \eta^{\symdex{2}{2}}_{\symdex{1}{2}} & \eta^{\symdex{2}{2}}_{\symdex{2}{2}}
    \end{matrix} \right)
  =\left( \begin{matrix}
      2\eta^{1}_{1} & 2 \eta^{1}_{2} & 0 \\
      \eta^{2}_{1} &  \eta^{1}_{1} + \eta^{2}_{2} & \eta^{1}_{2} \\
      0 & 2 \eta^{2}_{1} & 2\eta^{2}_{2}
    \end{matrix} \right) , \]
\[ \left(
    \begin{matrix}
      \Omega^{I}_{cd}
    \end{matrix} \right) =
  \left( \begin{matrix}
      \eta^{\tendex{1}{3}}_{\symdex{1}{1}} & \eta^{\tendex{1}{3}}_{\symdex{1}{2}} & \eta^{\tendex{1}{3}}_{\symdex{2}{2}}
      \\ \eta^{\tendex{1}{3'}}_{\symdex{1}{1}} & \eta^{\tendex{1}{3'}}_{\symdex{1}{2}} & \eta^{\tendex{1}{3'}}_{\symdex{2}{2}}
      \\ \eta^{\tendex{2}{3}}_{\symdex{1}{1}} & \eta^{\tendex{2}{3}}_{\symdex{1}{2}} & \eta^{\tendex{2}{3}}_{\symdex{2}{2}}
      \\ \eta^{\tendex{2}{3'}}_{\symdex{1}{1}} & \eta^{\tendex{2}{3'}}_{\symdex{1}{2}} & \eta^{\tendex{2}{3'}}_{\symdex{2}{2}}
    \end{matrix} \right)
  =\left( \begin{matrix}
      \eta_{\tendex{1}{3'}} & \eta_{\tendex{2}{3'}} & 0
      \\ \eta_{\tendex{1}{3}} & \eta_{\tendex{2}{3}} & 0
      \\ 0 &\eta_{\tendex{1}{3'}} & \eta_{\tendex{2}{3'}}
      \\ 0 &\eta_{\tendex{1}{3}} & \eta_{\tendex{2}{3}}
    \end{matrix} \right) , \]
and
\[ \left(
    \begin{matrix}
      \Omega^{I}_{J}
    \end{matrix} \right) =
  \left( \begin{matrix}
      \eta^{\tendex{1}{3}}_{\tendex{1}{3}} & \eta^{\tendex{1}{3}}_{\tendex{1}{3'}}
      & \eta^{\tendex{1}{3}}_{\tendex{2}{3}} & \eta^{\tendex{1}{3}}_{\tendex{2}{3'}}
      \\ \eta^{\tendex{1}{3'}}_{\tendex{1}{3}} & \eta^{\tendex{1}{3'}}_{\tendex{1}{3'}}
      & \eta^{\tendex{1}{3'}}_{\tendex{2}{3}} & \eta^{\tendex{1}{3'}}_{\tendex{2}{3'}}
      \\ \eta^{\tendex{2}{3}}_{\tendex{1}{3}} & \eta^{\tendex{2}{3}}_{\tendex{1}{3'}}
      & \eta^{\tendex{2}{3}}_{\tendex{2}{3}} & \eta^{\tendex{2}{3}}_{\tendex{2}{3'}}
      \\ \eta^{\tendex{2}{3'}}_{\tendex{1}{3}} & \eta^{\tendex{2}{3'}}_{\tendex{1}{3'}}
      & \eta^{\tendex{2}{3'}}_{\tendex{2}{3}} & \eta^{\tendex{2}{3'}}_{\tendex{2}{3'}}
    \end{matrix} \right)
  =\left( \begin{matrix}
      \eta^{1}_{1} + \eta^{3}_{3} & \eta^{3}_{3'}
      & \eta^{1}_{2} & 0
      \\ \eta^{3'}_{3} & \eta^{1}_{1} - \eta^{3}_{3}
      & 0 & \eta^{1}_{2}
      \\ \eta^{2}_{1} & 0
      & \eta^{2}_{2} + \eta^{3}_{3} & \eta^{3}_{3'}
      \\ 0 & \eta^{2}_{1}
      & \eta^{3'}_{3} & \eta^{2}_{2} - \eta^{3}_{3}
    \end{matrix} \right) , \]
and the torsion terms are given by
\begin{equation*}
    \left( \begin{matrix}
    T^{\vartheta^{\symdex{a}{b}}} \\ T^{\varpi^{\tendex{a}{i}}} \\ T^{\varpi^{\tendex{a}{i'}}}
    \end{matrix} \right) =
    \left( \begin{matrix}
    \varpi^{\tendex{a}{3}}\w\varpi^{\tendex{b}{3'}} + \varpi^{\tendex{b}{3}}\w\varpi^{\tendex{a}{3'}} \\
    0 \\
    0
    \end{matrix} \right) ,
\end{equation*}
so that the torsion functions are
\begin{equation}\label{eq: N torsion functions}
   T^{11}_{\tendex{1}{3},\tendex{1}{3'}} = 2, \quad
T^{12}_{\tendex{1}{3},\tendex{2}{3'}} = 1, \quad
T^{12}_{\tendex{2}{3},\tendex{1}{3'}} = 1, \quad
T^{22}_{\tendex{2}{3},\tendex{2}{3'}} = 2 ,
\end{equation}
skew in the lower indices.

From the structure equations it is clear that the forms \( \vartheta^{ab} \) define a corank 3 distribution on \( G_{N} \) that is invariant under the action of \( G^{0}_{N} \), so descends to a rank \( 4 \) distribution \( D' \) on \( N \).
Since the distribution \( D' \) is invariant, it is of constant symbol type.
In fact, one could go the other way (as was done for \( (2,3,5) \)-manifolds) and show that \( N \) is the model for the flat \( (4,7) \)-distribution of this symbol type.

There is a geometric description of \( D' \) that can be made directly on \( N \).
To do so, fix an isotropic \( 2 \)-plane \( E \in N \).
The \( \alpha \)-complement \( E^{\perp} \) of \( E \) is \( 4 \)-dimensional and contains \( E \).
By definition, \( \alpha \) descends to a well defined dual pairing \( E \otimes (\R^{6} / E^{\perp}) \to \R \), so that \( E^{\vee} \cong \R^{6} / E^{\perp} \) canonically.
As well, \( \alpha \) descends to a well defined symplectic form on \( Q := E^{\perp} / E \).

As \( N \) is a submanifold of the Grassmannian of \( 2 \)-planes in \( \R^{6} \), the tangent space \( T_{E} N \) is naturally a subspace of \( \Hom(E, \R^{6} / E) \), and is characterized as the kernel of the composition
\[ \Hom\left(E, \R^{6} / E\right) \to \Hom\left(E, \R^{6} / E^{\perp}\right) \cong E^{\vee} \otimes E^{\vee} \to \Wedge^{2} E^{\vee} . \]
There is a natural \( 4 \)-dimensional subspace of \( T_{E} N \), namely the subspace \( \Hom(E, E^{\perp} / E) \cong E^{\vee} \otimes Q \), and this choice at each point of \( N \) defines the invariant distribution \( D' \) on \( N \).
Note that the quotient \( T_{E} N / D'_{E} \) is isomorphic to \( \Sym^{2} E^{\vee} \).

From this, the vector space structure of \( \frg_{{N,-}} \) is clear.
Letting \( U' = E^{\vee} \), make the identification
\[ \frg_{N,-} = \frg_{N, -2} \oplus \frg_{N,-1} := \Sym^{2} U' \oplus \left(U'\otimes Q\right) . \]
This geometric description allows to fix once and for all a basis on \( \frg_{N,-} \).
Supposing without loss of generality that \( E \) is the isotropic plane spanned by the vectors \( e^{1}, e^{2} \) in \( \R^{6} \),
the basis \eqref{eq: symplectic basis} induces respectively bases \( e_{1}, e_{2} \) on \( U' \), then \( e_{3}, e_{3'} \) on \( Q \), and \( e^{1}, e^{2} \) on \( (U')^{\vee} \).
Using these, define the vectors
\begin{align*}
   e_{\tendex{a}{3}} = e_{a} \otimes e_{3},\quad e_{\tendex{a}{3'}} = e_{a} \otimes e_{3'}
   & \quad \mbox{ on } \quad U' \otimes Q \\
   e_{ab} = e_{a} \circ e_{b}
   & \quad \mbox{ on } \quad \Sym^{2} U'
\end{align*}
for \( a, b = 1, 2 \).
These define a basis
\begin{equation}\label{eq: nilpotent Iso2 model basis}
   e_{ab}, \quad e_{\tendex{a}{3}}, e_{\tendex{a}{3'}} 
\end{equation}
for \( \frg_{N,-} \), in agreement with the basis taken above.
Here and in the following, we use indices \( ab, cd, \ldots \) that range over symmetric pairs in \( \{1, 2\} \), and we will also employ the index set 
\[ I = \tendex{1}{3},\enspace \tendex{1}{3'},\enspace \tendex{2}{3},\enspace \tendex{2}{3'} . \]

Finally, it is not difficult to see in these bases that the Levi part \( G_{N,0} \) of \( G_{N} \) can be identified as \( \GL(U') \times \Sp(Q) \).

\section{Embeddable (2,3,5)-manifolds}\label{sec:Embeddable (2,3,5)-manifolds}
\subsection{The main theorem}\label{sec:The main theorem}
We can at this point state precise conditions for there to exist a Pfaffian embedding of \( (M, D) \) into \( (N, D') \),
where \( (M, D) \) is a \( (2,3,5) \)-manifold with \( D \) its derived rank \( 3 \) distribution,
and \( N \) is the isotropic Grassmannian with invariant \( 4 \)-distribution \( D' \).
We choose to do so in terms of the existence of a principal reduction of \( \B_{M} \),
where the structure group of the reductions will be the subgroup \( G^{(2)}_{M} \subset G^{0}_{M} \) with Lie algebra spanned by the vectors 
\[ X_{\zeta_{1}}, X_{\zeta_{2}}, X_{\gamma^{2}_{1}}, X_{\gamma^{0}_{1}} . \]
That these vectors are bracket-closed can be quickly checked by noting that the complementary forms
\[ \theta^{a}, \omega^{i}, \gamma^{1}_{2}, \gamma^{0}_{2}, \gamma, \gamma_{1}, \gamma_{2} \]
define a Frobenius system on \( G_{M} \).

\begin{theorem}\label{thm: main embeddability theorem}
The manifold \( (M, D) \) admits Pfaffian embeddings into \( (N, D') \) if and only if \( \B_{M} \) admits a \( G^{(2)}_{M} \)-principal reduction \( \mathcal{R} \subset \B_{M} \) such that the following relations hold,
 \begin{equation}\label{eq: main theorem reduction forms}
    \begin{aligned}
      \gamma^{1}_{2} & = 0, \\
      \gamma^{0}_{2} & = \tfrac{1}{14} A_{3} \theta^{1}, \\
      \gamma & = - \tfrac{4}{7} B_{3} \theta^{1} - \tfrac{5}{7} A_{3} \omega^{0}, \\
      \gamma_{2} & = - \tfrac{17}{7} C_{2} \theta^{1} + \tfrac{17}{14} A_{3} \omega^{1'}, \\
      \gamma_{1} & = -\tfrac{5}{7} C_{3} \theta^{1} + C_{2} \theta^{2} - \tfrac{22}{7} B_{3} \omega^{0} + \tfrac{9}{7} A_{4} \omega^{1'} + \tfrac{37}{14} A_{3} \omega^{2'}, 
    \end{aligned}
  \end{equation}
and the curvature functions on \( \mathcal{R} \) satisfy
\begin{align*}
  A_{4;1'} = {}& -5 B_{4} , \\
  A_{5;0;1'} = {}& 21 A_{5;1} .
\end{align*}
\end{theorem}

For such a reduction to exist requires that most of the curvatures vanish, as can be easily checked by comparing the differentials on the left and right side of each equation in display \eqref{eq: main theorem reduction forms}.
For example, one has
\begin{align*}
  \d \gamma^{1}_{2} = &
  \begin{multlined}[t]
  -\zeta_{2}\w\gamma^{1}_{2} {} + \gamma_{2}\w\theta^{1} + 3\gamma^{0}_{2}\w\omega^{1'}
  + A_{1} \theta^{2}\w\omega^{2'}
  + A_{2}(\theta^{1}\w\omega^{2'} + \theta^{2}\w\omega^{1'})
  \\
  + A_{3} \theta^{1}\w\omega^{1'}
  - 2 B_{1} \theta^{2}\w\omega^{0}
  - 2 B_{2} \theta^{1}\w\omega^{0} + C_{1} \theta^{1}\w\theta^{2}
  \end{multlined}
  \\ = & 
A_{1} \theta^{2} \w \omega^{2'}
+ A_{2} (\theta^{1} \w \omega^{2'} + \theta^{2} \w \omega^{1'})
- 2 B_{1} \theta^{2} \w \omega^{0} 
- 2 B_{2} \theta^{1} \w \omega^{0} 
+ C_{1} \theta^{1} \w \theta^{2} ,
\end{align*}
the first equality in general, and the second after restriction to \( \mathcal{R} \).
From this it is clear that \( \mathcal{R} \) is contained in the simultaneous zero locus of \( A_{1}, A_{2}, B_{1}, B_{2}, \) and \( C_{1} \).
In particular, for a generic \( M \), the common zero set is empty, and no such reduction can be made.

Similar constraints follow from the other equations, including higher order consequences.
In the interest of space, we list here only the first order ones, which are
\begin{equation}\label{eq: curvature consequences of reduction}
\begin{gathered}
   A_{1} = A_{2} = B_{1} = B_{2} = C_{1} = D_{1} = \tilde{D}_{1} = \tilde{D}_{2} = \tilde{E}_{1} = \tilde{F}_{1} = 0, \\
   E = \tilde{E}_{2} = \tfrac{9}{14} (A_{3})^2,
   \quad \tilde{D}_{3} = -\tfrac{2}{3} D_{2} .
\end{gathered}
\end{equation}
The higher order curvature consequences can be calculated as needed.

\subsection{The coframed embedding PDE}\label{sec:The coframed embedding PDE}
To prove Theorem \ref{thm: main embeddability theorem}, it suffices to determine when the exterior differential system \( (P, \I) \) admits solutions, where we recall that \( (P, \I) \) was defined in Subsection \ref{sec:The 1-jet contact system and the embedding PDE}.
In this analysis, it will be profitable to include the coframing data that the Cartan geometries on \( M \) and \( N \) allow.
To this end, we describe an enlarged system \( (\BP, \I) \) that admits integral sections if and only if \( P \) does.

Let \( (\mathcal{B}, \Omega) \) be the product Cartan geometry on \( M \times N \), so that \( \mathcal{B} \) is a \( G^{0} := G^{0}_{M} \times G^{0}_{N} \)-principal bundle over \( M \times N \) and \( \Omega \) the induced \( \frg := \frg_{M} \oplus \frg_{N} \)-valued Cartan connection.
With this, define the space of \emph{coframed jets} \( \mathcal{B} J^{1} \) over \( M \times N \) by the following pullback diagram.
\begin{equation}\label{eq: jet coframe pullback definition}
  \begin{tikzcd}
    \mathcal{B} J^{1} \ar[r] \ar[d] & \mathcal{B} \ar[d] \\
    J^{1} \ar[r] & M \times N  
  \end{tikzcd}
\end{equation}
On \( \BJets \), denote by the same symbols the pullbacks of both the Cartan form \( \Omega \) and the contact form \( \Theta \), as well as the differentially closed ideal \( \I \) generated by the component \( 1 \)-forms of \( \Theta \).

We claim that locally, integral sections of \( \BJets \) exist if and only if they exist for \( J^{1} \).
It is clear that any integral section of \( \mathcal{B} J^{1} \to M \) pushes down to an integral section of \( J^{1} \to M \), simply by composing with the projection \( \mathcal{B} J^{1} \to J^{1} \).
Conversely, since there always exist local sections of \( \mathcal{B} \to M \times N \), and any such section pulls back to a section of \( \mathcal{B} J^{1} \to J^{1} \), any integral section of \( J^{1} \) can locally be lifted to an integral section of \( \mathcal{B} J^{1} \).
We remark in particular, that an integral section of \( \BJets \to M \) is essentially the data of an integral section of \( J^{1} \) plus an arbitrary choice of coframing data.

The claim can also be seen by noting that the fibers of \( \BJets \to J^{1} \) are simply the Cauchy leaves of the differential system \( \I \) on \( \BJets \).
This being so, the integral sections of \( (J^{1}, \I) \) (equivalently: the integral submanifolds that submerse to \( M \)) are in local bijection with the maximal integral submanifolds of \( (\BJets, \I) \) that submerse to \( M \).
For given the former, its preimage in \( \BJets \) will be integral, and given the latter, its projection to \( J^{1} \) will be integral.

From either perspective it can be seen that the claim will hold even after restriction to submanifolds of \( \BJets \) transverse to the Cauchy directions.
This is one advantage of working on \( \BJets \), and will allow to make adaptations of coframing data to solutions at the \( 1 \)-jet level, simply by considering submanifolds of \( \BJets \).

A second advantage is that one can use the coframing form \( \Omega \) on \( \BJets \) to improve the bundle-valued jet linearization map \eqref{eq: jet linearization map} to a function taking values in a single vector space \( \Hom(\frg_{M,-}, \frg_{N,-}) \).
Indeed, points of \( \BJets \) are the same as the triples \( (j^{1}_{x}(f), p, p') \in J^{1} \times \B_{M} \times \B_{N} \) for which \( p \) is in the fiber of \( x \) and \( p' \) is in the fiber of \( f(x) \).
Recalling that \( p \) determines an isomorphism \( \varphi_{M,p} \colon T_{x} M \to \frg_{M,-} \) and \( p' \) an isomorphism \( \varphi_{N,p'} \colon T_{f(x)} N \to \frg_{N,-} \), define the \emph{canonical jet linearization map},
\[ \begin{tikzcd}[row sep={8mm,between origins}]
  \jhom \colon \BJets \ar[r] & \Hom(\frg_{M,-}, \frg_{N,-}) \\
  (j^{1}_{x}(f), p, p') \ar[r, mapsto] & \varphi_{N,p'} \circ \d f_{x} \circ \varphi^{-1}_{M,p} .
\end{tikzcd} \]

There are two natural group actions on \( \BJets \), and \( \jhom \) is compatible with both.
The first action is the right principal action of \( G^{0} = G^{0}_{M} \times G^{0}_{N} \) on the fibers of \( \BJets \to J^{1} \).
It is not difficult to check that \( H \) is equivariant for this action if we fix the tensor product right action of \( G^{0} \) on \( \Hom(\frg_{M,-}, \frg_{N,-}) \), so that
\[ \jhom\left((j^{1}_{x}(f), p, p') \cdot (g,h) \right)
= \Ad_{h^{-1}} \circ \varphi_{N,p'} \circ \d f \circ \varphi_{M,p}^{-1} \circ \Ad_{g}
= \jhom\left((j^{1}_{x}(f), p, p')\right)\cdot (g,h) . \]
The second action on \( \BJets \) is the fiber product of the left actions of \( G_{N} \) on \( J^{1} \) and \( \B \), where the action on \( \B \) is clear and the action on \( J^{1}(M, N) \) is by postcomposition with left multiplication, \( g \cdot j^{1}_{x}(f) = j^{1}_{x}(L_{g} \circ f) \).
It is again not difficult to check that \( H \) is invariant for this action,
\[ H\left(g \cdot (j^{1}(f), p, p')\right)
= \varphi_{N,gp'} \circ \d L_{g} \circ \d f \circ \varphi_{M, p}^{-1}
= \varphi_{N,p'} \circ \d f \circ \varphi_{M, p}^{-1}
= H\left((j^{1}(f), p, p')\right) . \]

Returning to the embedding PDE, restrict the coframed jet bundle over the space \( P \subset J^{1} \), to define the space \( \mathcal{B} P \), as in the following diagram.
\[  \begin{tikzcd}
    \mathcal{B} P \ar[r,hook]\ar[d] & \mathcal{B} J^{1} \ar[d] \\
    P \ar[r, hook] & J^{1} 
  \end{tikzcd} \]
The point then is that locally \emph{there exist integral sections of \( P \to M \) if and only if there exist integral sections of \( \mathcal{B} P \to M \).}
(And we remind, the existence of integral sections of \( P \) is the same as the existence of Pfaffian embeddings \( M \to N \).)

Upon restriction to \( \BP \), the canonical jet linearization map takes values in
\begin{equation}\label{eq: first tableau}  
\left(
  \begin{matrix}
    \Sym^2 U' \otimes U^\vee & 0 \\
    (U' \otimes Q) \otimes U^\vee &
    (U' \otimes Q) \otimes (\R \oplus U)^\vee
  \end{matrix} \right) \subset \Hom(\frg_{M,-}, \frg_{N,-}) .
\end{equation}
Relative the bases of \( \frg_{M,-} \) and \( \frg_{N,-} \) fixed in Equations \eqref{eq: nilpotent 235 model basis} and \eqref{eq: nilpotent Iso2 model basis}, the block components of \( \jhom \) decompose as
\begin{align*}
  (\tablA{a}{b}{c})
={}& \left( \begin{matrix}
    \tablA{1}{1}{1} & \tablA{1}{1}{2} \\
    \tablA{1}{2}{1} & \tablA{1}{2}{2} \\
    \tablA{2}{2}{1} & \tablA{2}{2}{2} \end{matrix} \right), \\
(\tablBSingle{I}{c})
={}& \left( \begin{matrix}
    \tablB{1}{3 }{1} & \tablB{1}{3 }{2} \\
    \tablB{1}{3'}{1} & \tablB{1}{3'}{2} \\
    \tablB{2}{3 }{1} & \tablB{2}{3 }{2} \\
    \tablB{2}{3'}{1} & \tablB{2}{3'}{2}
  \end{matrix} \right), \qquad
(\tablCSingle{I}{i})
= \left( \begin{matrix}
    \tablC{1}{3 }{0} & \tablC{1}{3 }{1'} & \tablC{1}{3 }{2'} \\
    \tablC{1}{3'}{0} & \tablC{1}{3'}{1'} & \tablC{1}{3'}{2'} \\
    \tablC{2}{3 }{0} & \tablC{2}{3 }{1'} & \tablC{2}{3 }{2'} \\
    \tablC{2}{3'}{0} & \tablC{2}{3'}{1'} & \tablC{2}{3'}{2'}
\end{matrix} \right) . 
\end{align*}

In terms of the jet linearization map \( \jhom \), the submanifold \( \BP \) is cut out by the system of equations \( \jhom^{\symdex{a}{b}}_{i} = 0 \).
This system of equations is manifestly invariant for both the right \( G^{0} \)-action and the left \( G_{N} \)-action.
More generally, one may restrict to solutions of equations in the components of \( \jhom \) that are not \( G^{0} \)-invariant.
In doing so, one has restrictions transverse to the Cauchy directions, as discussed previously in this subsection, and so may obtain normalizations at the cost of reducing the action of \( G^{0} \) to the isotropy group of solutions.
See Remark \ref{rmk: V1 normalization} for the first example of this.
On the other hand, any such submanifold is still the union of orbits of \( G_{N} \).

The contact form \( \Theta \) on \( \BJets \) is easily expressed in terms of \( H \) and \( \Omega \).
For each point \( p \in \BJets \), the corresponding isomorphism of \( T_{y} N \) with \( \frg_{N,-} \) allows one to regard \( \Theta \) as a \( \frg_{N,-} \)-valued form.
Relative the fixed basis of \( \frg_{N,-} \), the component \( 1 \)-forms of \( \Theta \) may be denoted by \( \sTheta{a}{b}, \tThetaSingle{I} \).
Then it follows from the definitions that
\[ \sTheta{a}{b} = \vartheta^{\symdex{a}{b}} - \tablA{a}{b}{c} \theta^{c} - \jhom^{\symdex{a}{b}}_{i} \omega^{i} , \qquad
   \tThetaSingle{I} = \varpi^{I} - \tablBSingle{I}{a} \theta^{a} - \tablCSingle{I}{i}\omega^{i} , \]
and that upon restriction to \( \BP \) these simplify to
\begin{equation}\label{eq: first framed contact forms}
   \sTheta{a}{b} = \vartheta^{\symdex{a}{b}} - \tablA{a}{b}{c} \theta^{c} , \qquad
   \tThetaSingle{I} = \varpi^{I} - \tablBSingle{I}{a} \theta^{a} - \tablCSingle{I}{i}\omega^{i} . 
\end{equation}

  The differential of \( \Theta \) on \( \BP \) is now straightforward to compute,\hspace{-0.7mm}\footnote{In fact, the following structure equations hold more generally, and are simpler to compute, for general coframed jet bundles.
  If \( M \) has coframing and connection forms \( \Omega^{i}, \Omega^{i}_{j} \) and \( N \) has coframing and connection forms \( \Omega^{I}, \Omega^{I}_{J} \), with respective torsions \( T^{i}_{j,k} \) and \( T^{I}_{J,K} \), then one finds that
  \( \pi^{I}_{i} = \d \jhom^{I}_{i} + \Omega^{I}_{J} \jhom^{J}_{i} - \jhom^{I}_{j} \Omega^{j}_{i} \) and \( T^{I}_{i,j} = T^{I}_{J,K} \jhom^{J}_{i} \jhom^{K}_{j} - \jhom^{I}_{k} T^{k}_{i,j} \).
  The specific structure equations here are just a specialization of these equations.} and one finds
  
  \begin{equation}\label{eq: Theta structure equations}
    \begin{aligned}
      \d\sTheta{a}{b} {}&
      \equiv - \ftablA{a}{b}{e} \w \theta^{e}
      + T^{\symdex{a}{b}}_{1,2}\ \theta^{1}\w\theta^{2}
      + T^{\symdex{a}{b}}_{e,i}\ \theta^{e}\w\omega^{i}
      + T^{\symdex{a}{b}}_{i,j}\ \omega^{i}\w\omega^{j}
      & \mod{\Theta^{cd},\Theta^{J}} \\
      \d\tThetaSingle{I} {}&
      \equiv - \ftablBSingle{I}{a}\w \theta^{a} - \ftablCSingle{I}{i}\w \omega^{i}
      + T^{\tendexSingle{I}}_{i,j} \omega^{i}\w\omega^{j}
      & \mod{\Theta^{cd},\Theta^{J}} ,
    \end{aligned}
  \end{equation}
  with \emph{tableau forms} \( \pi \) given by
  \begin{equation}\label{eq: tableau forms definition}
    \begin{aligned}
      \ftablA{a}{b}{c} &= \d \tablA{a}{b}{c} + \Omega^{\symdex{a}{b}}_{\symdex{d}{e}} \tablA{d}{e}{c} - \tablA{a}{b}{d} \Omega^{\thindex{d}}_{\thindex{c}}
      \\ \ftablBSingle{I}{a}
      &= \d \tablBSingle{I}{a}
      +\Omega^{\tendexSingle{I}}_{\symdex{b}{c}} \tablA{b}{c}{a}
      +\Omega^{\tendexSingle{I}}_{\tendexSingle{J}} \tablBSingle{J}{a}
      -\tablBSingle{I}{c} \Omega^{\thindex{c}}_{\thindex{a}} - \tablCSingle{I}{0} \Omega^{0}_{\thindex{a}} - \tablCSingle{I}{b'} \Omega^{b'}_{\thindex{a}}
      \\ \ftablCSingle{I}{0}
      &= \d \tablCSingle{I}{0}
      +\Omega^{\tendexSingle{I}}_{\tendexSingle{J}} \tablCSingle{J}{0}
      -\tablCSingle{I}{0} \Omega^{0}_{0}
      -\tablCSingle{I}{b'}\Omega^{b'}_{0}
      \\ \ftablCSingle{I}{a'}
      &= \d \tablCSingle{I}{a'}
      +\Omega^{\tendexSingle{I}}_{\tendexSingle{J}} \tablCSingle{J}{a'}
      -\tablCSingle{I}{b'}\Omega^{b'}_{a'} ,
    \end{aligned}
    \end{equation}
  and mixed torsion functions \( T \) given by the combinations of torsion functions defined in Equations \eqref{eq: M torsion functions} and \eqref{eq: N torsion functions},

  \begin{align*}
     T^{\symdex{a}{b}}_{1,2} & = T^{\symdex{a}{b}}_{\tendex{c}{3}, \tendex{d}{3'}}\tablB{c}{3}{1}\tablB{d}{3'}{2}, \\
     T^{\symdex{a}{b}}_{e,j} & = T^{\symdex{a}{b}}_{\tendex{c}{3}, \tendex{d}{3'}}\tablB{c}{3}{e}\tablB{d}{3'}{j}, \\
     T^{\symdex{a}{b}}_{i,j} & = T^{\symdex{a}{b}}_{\tendex{c}{3}, \tendex{d}{3'}}\tablB{c}{3}{i}\tablB{d}{3'}{j} - \tablA{a}{b}{c} T^{c}_{i,j} , \\
     T^{I}_{i,j} & = - \tablBSingle{I}{c} T^{c}_{i,j} .
  \end{align*}

  The torsion functions can be usefully re-expressed in terms of the symplectic form \( \alpha \) of \( Q \) if one considers the \( (U' \otimes Q) \otimes \frg_{M,-}^\vee \)-valued component of the jet linearization map \( \jhom \) as a collection of \( Q \)-valued functions,
  \begin{equation}\label{eq: Q-valued jet linearization}
  q^{a}_{b} = \left(
  \begin{matrix}
  \tablC{a}{3}{b} \\ \tablC{a}{3'}{b}
  \end{matrix} \right),
  \qquad
  q^{a}_{i} = \left(
  \begin{matrix}
  \tablC{a}{3}{i} \\ \tablC{a}{3'}{i}
  \end{matrix} \right) , 
\end{equation}
  so that
  \[ (\tablBSingle{I}{a}, \tablCSingle{I}{i})
    = \left( \begin{matrix}
        q^{1}_{1} & q^{1}_{2} & q^{1}_{0} & q^{1}_{1'} & q^{1}_{2'} \\
        q^{2}_{1} & q^{2}_{2} & q^{2}_{0} & q^{2}_{1'} & q^{2}_{2'} \\
      \end{matrix} \right) . \]
  In this notation, the torsion functions can be written more explicitly as
  \begin{align*}
     T^{\symdex{a}{b}}_{1,2} & = \alpha(q^{a}_{1}, q^{b}_{2}) + \alpha(q^{b}_{1}, q^{a}_{2}) , &
     T^{\symdex{a}{b}}_{0,c'} & = \alpha(q^{a}_{0}, q^{b}_{c'}) + \alpha(q^{b}_{0}, q^{a}_{c'}) - 3\tablA{a}{b}{c} , \\
     T^{\symdex{a}{b}}_{c,j} & = \alpha(q^{a}_{c}, q^{b}_{j}) + \alpha(q^{b}_{c}, q^{a}_{j}) , &
     T^{\symdex{a}{b}}_{1',2'} & = \alpha(q^{a}_{1'}, q^{b}_{2'}) + \alpha(q^{b}_{1'}, q^{a}_{2'}) , \\
     T^{I}_{0,c'} & = - 3 \tablBSingle{I}{c} , &
     T^{I}_{1',2'} & = - 2 \tablCSingle{I}{0} .
  \end{align*}

\subsection{The first integrability condition}\label{sec:The first integrability condition}
It follows immediately from the structure Equations \eqref{eq: Theta structure equations} that \( (\BP, \I) \) is a \emph{linear Pfaffian system}, so that we may apply the full theory as in \cite{Bryant.Chern.eaExteriorDifferentialSystems1991}, chapter IV.
In particular, the forms
\begin{equation}\label{eq: BP partial Pfaffian coframing}
    \sTheta{a}{b}, \tThetaSingle{I}, \quad \theta^{a}, \omega^{i}, \quad
    \ftablA{a}{b}{c}, \ftablBSingle{I}{c} ,\ftablCSingle{I}{i}
\end{equation}
define a \emph{partial linear Pfaffian coframing}, with the forms
\[ \Theta^{ab}, \Theta^{I}, \theta^{a}, \omega^{i} \]
spanning an \emph{independence condition} for the linear Pfaffian system.
Any integral section of \( (\BP, \I) \) will by definition pull back \( \Theta^{ab}, \Theta^{I} \) to zero and pull back \( \theta^{a}, \omega^{i} \) to a coframing of \( M \).
We remark that the linear Pfaffian coframing is \emph{partial} only because it doesn't span, reflecting the fact that the Pfaffian system has Cauchy characteristic directions.
These are useful to carry along, as they represent symmetries of the system, and otherwise cause no harm, as discussed after the definition of \( \BJets \).

The tableau forms \( \pi \) are \( \Hom(\frg_{M,-}, \frg_{N,-}) \)-valued, and the \emph{tableau bundle} of \( (\BP, \I) \) can be identified as a subbundle of the trivial \( \Hom(\frg_{M,-}, \frg_{N,-}) \)-bundle over \( \BP \).
Precisely, the fiber of the tableau bundle over \( p \in \BP \) is given by evaluating all the forms in display \eqref{eq: BP partial Pfaffian coframing} against all vectors \( X \in T_{p} \BP \) that are annihilated by the independence condition.
It is immediate then that all fibers of the tableau bundle are isomorphic, being of the block form in Equation \eqref{eq: first tableau}.
The \emph{symbol} of the linear Pfaffian system at each point, the subspace of \( \Hom(\frg_{M,-}, \frg_{N,-})^{\vee} \) annihilating the tableau, is thus easily seen to be spanned by the matrices \( e_{\symdex{a}{b}} \otimes e^{i} \) with respect to the bases fixed on \( \frg_{M,-} \) and \( \frg_{N,-} \) above.

This non-trivial symbol, or equivalently the fact that the block \( \pi^{\symdex{a}{b}}_{i} \) of the tableau forms vanishes, gives rise to the first potential \emph{integrability conditions}, obstructions to integral sections.
Indeed, it is immediate that 
\[ \d \sTheta{a}{b} \equiv T^{\symdex{a}{b}}_{i,j}\omega^{i}\w\omega^{j} \mod{\theta^{c}, \sTheta{c}{d}, \tThetaSingle{I}} , \]
so the standard argument applies, showing that integral sections of \( \BP \to M \) cannot pass through any point of \( \BP \) where the functions \( T^{\symdex{a}{b}}_{i,j} \) are non-zero.
We are thus obliged to restrict to the zero set of \( T^{\symdex{a}{b}}_{i,j} \), which is a \( 52 \)-dimensional manifold, to be denoted \( \V'_{1} \).

\section{First order reductions}\label{sec:First order reductions}
\subsection{Normalization after the first integrability condition}\label{sec:Normalization after the first integrability condition}
The integrability conditions described in the previous section comprise \( 9 \) equations in total, namely,
\begin{equation}\label{eq: first integrability equations}
    T^{\symdex{1}{1}}_{0,1'} = T^{\symdex{1}{1}}_{0,2'} = T^{\symdex{1}{1}}_{1',2'} =
  T^{\symdex{1}{2}}_{0,1'} = T^{\symdex{1}{2}}_{0,2'} = T^{\symdex{1}{2}}_{1',2'} =
  T^{\symdex{2}{2}}_{0,1'} = T^{\symdex{2}{2}}_{0,2'} = T^{\symdex{2}{2}}_{1',2'} = 0 , 
\end{equation}
establishing a system of fiber-wise equations in the coefficients \( \tablA{a}{b}{c} \) and \( \tablC{a}{3}{i}, \tablC{a}{3'}{i} \) of \( \jhom \).
Upon restriction to \( \V'_{1} \), the tableau bundle loses \( 9 \) dimensions, and correspondingly, the symbol gains \( 9 \) relations.
Before describing these new symbol relations, it will be useful to normalize, which is carried out in the following Proposition.

\begin{proposition}\label{thm: prop normalization of tableau}
  For each point \( p \in \V_{1}' \), there exists an element \( g \in G^{0} = G^{0}_{M} \times G^{0}_{N} \) so that the components of the jet linearization function evaluated at \( p \cdot g \) are given by
  \begin{equation}\label{eq: normalized A}
    \left( \begin{matrix}
        \tablA{1}{1}{1} & \tablA{1}{1}{2} \\
        \tablA{1}{2}{1} & \tablA{1}{2}{2} \\
        \tablA{2}{2}{1} & \tablA{2}{2}{2}
      \end{matrix} \right) = \left( \begin{matrix}
        0 & \tfrac{2}{3} \\ \tfrac{1}{3} & 0 \\ 0 & 0
      \end{matrix} \right), \hfill
  \end{equation}
  and
  \begin{equation}
    \label{eq: normalized B, C}
    \left( \begin{matrix}
        \tablB{1}{3}{1} & \tablB{1}{3}{2} \\
        \tablB{1}{3'}{1} & \tablB{1}{3'}{2} \\
        \tablB{2}{3}{1} & \tablB{2}{3}{2} \\
        \tablB{2}{3'}{1} & \tablB{2}{3'}{2}
      \end{matrix} \right)
    = \left( \begin{matrix}
        0 & \tablB{1}{3}{2} \\
        0 & 0 \\
        0 & \tablB{2}{3}{2} \\
        0 & \tablB{2}{3'}{2}
      \end{matrix} \right),
       \quad
    \left( \begin{matrix}
        \tablC{1}{3}{0} & \tablC{1}{3}{1'} & \tablC{1}{3}{2'} \\
        \tablC{1}{3'}{0} & \tablC{1}{3'}{1'} & \tablC{1}{3'}{2'} \\
        \tablC{2}{3}{0} & \tablC{2}{3}{1'} & \tablC{2}{3}{2'} \\
        \tablC{2}{3'}{0} & \tablC{2}{3'}{1'} & \tablC{2}{3'}{2'}
      \end{matrix} \right)
    = \left( \begin{matrix}
        1 & 0 & 0 \\
        0 & 0 & 1 \\
        0 & 0 & 0 \\
        0 & 1 & 0
      \end{matrix} \right)
  \end{equation}
\end{proposition}

\begin{remark}\label{rmk: V1 normalization}
  Let \( \V_{1} \) denote the \( 38 \)-dimensional submanifold of \( \V'_{1} \) on which the normalizations \eqref{eq: normalized A} and \eqref{eq: normalized B, C} hold.
  It should be noted that the normalizations made here do not introduce new symbol relations into the partial differential equation, as they merely cut down on the Cauchy directions.
  Rather, these normalizations have the effect of adapting framings on \( M \) and \( N \) to the \( 1 \)-jets of solutions.
  In particular, the structure group acting on \( \V_{1} \) is the subgroup \( G^{(1)} \subset G^{0} \) that stabilizes the elements of \( \Hom(\frg_{M,-}, \frg_{N,-}) \) as determined in the statement of the Proposition.
\end{remark}

\begin{proof}
In the notation of Equation \eqref{eq: Q-valued jet linearization}, the integrability conditions are
\[ 3 \tablA{a}{b}{c} = \alpha(q^{a}_{0}, q^{b}_{c'}) + \alpha(q^{b}_{0}, q^{a}_{c'}),
  \qquad 0 = \alpha(q^{a}_{1'}, q^{b}_{2'}) + \alpha(q^{b}_{1'}, q^{a}_{2'}) . \]
  It follows from the second set of equations, with \( a = 1, b = 1 \), that the vectors \( q^{1}_{1'}, q^{1}_{2'} \) are linearly dependent, and likewise with \( a = 2, b = 2 \) that \( q^{2}_{1'}, q^{2}_{2'} \) are linearly dependent, so for any point \( p \in \V'_{1} \), we may write
\[ \tablCSingle{I}{i}(p) =
  \left( \begin{matrix}
      q_{1} & q_{3} \otimes v^{1} \\
      q_{2} & q_{4} \otimes v^{2}
  \end{matrix} \right) \]
for some vectors \( q_{1},q_{2},q_{3},q_{4} \in Q \) and \( v^{1}, v^{2} \in U^\vee \).
The \( a = 1, b = 2 \) component of the second set of equations now reads as \( 0 = \alpha(q_{3}, q_{4}) \det(v^{1}, v^{2}) \), and it follows that either \( v^{1}, v^{2} \) are linearly dependent or \( q_{3}, q_{4} \) are linearly dependent.
The remaining \( 6 \) integrability conditions require that
\[ \left( \begin{matrix}
      \tablA{1}{1}{1} & \tablA{1}{1}{2} \\
      \tablA{1}{2}{1} & \tablA{1}{2}{2} \\
      \tablA{2}{2}{1} & \tablA{2}{2}{2}
    \end{matrix} \right)(p)
  = \frac{1}{3}\left( \begin{matrix}
      2\alpha(q_{1}, q_{3}) v^{1} \\
      \alpha(q_{2}, q_{3}) v^{1} + \alpha(q_{1}, q_{4}) v^{2} \\
      2\alpha(q_{2}, q_{4}) v^{2}
\end{matrix} \right) .
\]
It follows that \( \jhom \) has rank less than \( 5 \) if either the \( v^{1}, v^{2} \) are linearly dependent or one of \( q_{3}, q_{4} \) vanishes.
This would contradict the assumption that \( \jhom \) is injective, so suppose \( q_{3}, q_{4} \) dependent but both non-vanishing, and then rescaling \( v^{1}, v^{2} \) as necessary, suppose further that \( q_{3} = q_{4} \).

Now, the action of \( G^{0}_{M} \times G^{0}_{N} \) on \( \tablCSingle{I}{i} \) factors through the group \( G_{M,0} \times G_{N,0} = \GL(U) \times \Sp(Q) \times \GL(U') \).
It is not difficult to use this action to arrange that
\[ q_{1} = \left(
    \begin{matrix}
      1 \\ 0
    \end{matrix} \right)
  \qquad \mbox{ and } \qquad
  q_{3} = \left(
    \begin{matrix}
      0 \\ 1
    \end{matrix} \right) , \]
then that \( q_{2} \) is a multiple of \( q_{3} \), and finally that
\[ v^{1} = (0, 1), \quad v^{2} = (1, 0) . \]
This exhausts the reductions that can be made using the \( G_{M,0} \) and \( G_{N,0} \) actions, demonstrating the existence of \( g \) such that Equation \eqref{eq: normalized A} holds, and 
\[  \left( \begin{matrix}
        \tablC{1}{3}{0} & \tablC{1}{3}{1'} & \tablC{1}{3}{2'} \\
        \tablC{1}{3'}{0} & \tablC{1}{3'}{1'} & \tablC{1}{3'}{2'} \\
        \tablC{2}{3}{0} & \tablC{2}{3}{1'} & \tablC{2}{3}{2'} \\
        \tablC{2}{3'}{0} & \tablC{2}{3'}{1'} & \tablC{2}{3'}{2'}
      \end{matrix} \right)
    (p \cdot g)
    = \left( \begin{matrix}
        1 & 0 & 0 \\
        0 & 0 & 1 \\
        0 & 0 & 0 \\
        \tablC{2}{3'}{0} & 1 & 0
      \end{matrix} \right) \]
with \( \tablC{2}{3'}{0} \) and all of the \( \tablBSingle{I}{a} \) as yet unnormalized.

The higher order components of \( G^{0} \) can be used to further reduce \( \jhom \).
First, on \( M \), the action by \( \exp(\tfrac{1}{2}\tablC{2}{3'}{0} X_{\gamma^{0}_{2}}) \) will preserve the existing normalizations and reduce \( \tablC{2}{3'}{0} \) to zero.
Similarly for the action of \( \exp(\tablB{1}{3'}{2} X_{\gamma}) \), which reduces \( \tablB{1}{3'}{2} \) to zero.
Finally, using the transformations on \( N \), an action by the exponential of the element
\[
-3\tablB{2}{3'}{1} X_{\eta_{\tendex{1}{3}}} - 3\tablB{2}{3}{1} X_{\eta_{\tendex{1}{3'}}} - 3\tablB{1}{3'}{1} X_{\eta_{\tendex{2}{3}}} - 3\tablB{1}{3}{1} X_{\eta_{\tendex{2}{3'}}}
\]
of \( \frg_{N,1} \) ensures that all of the \( \tablBSingle{I}{1} \) vanish.
\end{proof}

\subsection{The tableau forms}\label{sec:The tableau forms}
After the normalizations of Proposition \ref{thm: prop normalization of tableau}, the jet tableau forms on \( \V_{1} \) simplify.
Explicitly, following Equations \eqref{eq: tableau forms definition}, it holds that
\[
\left( \begin{matrix}
\ftablA{1}{1}{1} & \ftablA{1}{1}{2} \\
\ftablA{1}{2}{1} & \ftablA{1}{2}{2} \\
\ftablA{2}{2}{1} & \ftablA{2}{2}{2}
\end{matrix} \right)
= \left( \begin{matrix}
    \tfrac{2}{3}\eta^{1}_{2} - \tfrac{2}{3}\gamma^{2}_{1} &
    \tfrac{4}{3}\eta^{1}_{1} - 2\zeta_{1} - \tfrac{2}{3}\zeta_{2} \\
    \tfrac{1}{3}\eta^{1}_{1} + \tfrac{1}{3}\eta^{2}_{2} - \zeta_{1} - \tfrac{2}{3}\zeta_{2} &
    \tfrac{2}{3}\eta^{2}_{1} - \tfrac{1}{3}\gamma^{1}_{2} \\
    \tfrac{2}{3}\eta^{2}_{1} &
    0
\end{matrix} \right)
\]
and
\begin{align*}
  \left(
    \begin{matrix}
      \ftablC{1}{3}{0} & \ftablC{1}{3}{1'} & \ftablC{1}{3}{2'}
      \\ \ftablC{1}{3'}{0} & \ftablC{1}{3'}{1'} & \ftablC{1}{3'}{2'}
      \\ \ftablC{2}{3}{0} & \ftablC{2}{3}{1'} & \ftablC{2}{3}{2'}
      \\ \ftablC{2}{3'}{0} & \ftablC{2}{3'}{1'} & \ftablC{2}{3'}{2'}
    \end{matrix} \right)
  = \left(
    \begin{matrix}
      \eta^{1 }_{1} + \eta^{3 }_{3} -2\zeta_{1} - \zeta_{2} & 0 & \eta^{3}_{3'}
      \\ \eta^{3'}_{3} + 2\gamma^{0}_{1}
      & \eta^{1}_{2} - \gamma^{2}_{1}
      & \eta^{1}_{1} - \eta^{3}_{3} - \zeta_{1}
      \\ \eta^{2}_{1} & \eta^{3}_{3'} & 0
      \\ - 2\gamma^{0}_{2}
      &  \eta^{2}_{2} - \eta^{3}_{3} - \zeta_{1} - \zeta_{2}
      &  \eta^{2}_{1} - \gamma^{1}_{2}
    \end{matrix} \right)
\end{align*}

The remaining component \( (\ftablBSingle{I}{a}) \) of the tableau forms is as yet messy to write.
But for sake of recording in one place, we note that further integrability conditions found below will require to restrict to a submanifold \( \V_{4} \) such that
\begin{equation}\label{eq: V4 equations}
   \tablB{1}{3}{2} = \tablB{2}{3}{2} = \tablB{2}{3'}{2} = 0 . 
\end{equation}
On this manifold the situation is cleaner, so that
\begin{equation*}
     \left(
    \begin{matrix}
      \ftablB{1}{3 }{1} & \ftablB{1}{3 }{2} \\
      \ftablB{1}{3'}{1} & \ftablB{1}{3'}{2} \\
      \ftablB{2}{3 }{1} & \ftablB{2}{3 }{2} \\
      \ftablB{2}{3'}{1} & \ftablB{2}{3'}{2}
    \end{matrix} \right)
  = \left(
    \begin{matrix}
      \tfrac{1}{3}\eta_{\tendex{2}{3'}} - \gamma^{0}_{1}
      &  \tfrac{2}{3}\eta_{\tendex{1}{3'}} - \gamma^{0}_{2}
      \\ \tfrac{1}{3}\eta_{\tendex{2}{3}}
      & \tfrac{2}{3}\eta_{\tendex{1}{3}} - \gamma
      \\ \tfrac{1}{3}\eta_{\tendex{1}{3'}}
      & 0
      \\ \tfrac{1}{3}\eta_{\tendex{1}{3}} - \gamma
      & 0
    \end{matrix} \right) .
\end{equation*}

\subsection{The new symbol relations}\label{sec:The new symbol relations}
Restriction by the nine integrability Equations \eqref{eq: first integrability equations} introduces nine new symbol relations between the tableau forms of the linear Pfaffian system \( (\V_{1}, \I) \).
This can be computed by taking the differentials of the integrability equations, keeping in mind the normalizations made in Proposition \ref{thm: prop normalization of tableau}.
The relations can also be checked directly from the manner in which the coframing forms of \( \B \) embed into the tableau, as described in the previous section.
The new relations are as follows,
\begin{equation}\label{eq: first integrability symbol relations}
  \begin{aligned}
    3\ftablA{1}{1}{1} & = 2\ftablC{1}{3'}{1'}
    & \qquad 3\ftablA{1}{2}{1} & = \ftablC{2}{3'}{1'} + \ftablC{1}{3}{0}
    & \qquad 3\ftablA{2}{2}{1} & = 2\ftablC{2}{3}{0} \\
    3\ftablA{1}{1}{2} & = 2\ftablC{1}{3}{0} + 2\ftablC{1}{3'}{2'}
    &  3\ftablA{1}{2}{2} & = \ftablC{2}{3'}{2'} + \ftablC{2}{3}{0}
    &  3\ftablA{2}{2}{2} & = 0 \\
    0 & = \ftablC{1}{3}{1'}
    &  0 & = \ftablC{1}{3}{2'} - \ftablC{2}{3}{1'}
    &  0 & = \ftablC{2}{3}{2'} .
  \end{aligned}
\end{equation}

Again for the sake of recording in one place, we note that on the submanifold \( \V_{4} \) defined by Equations \eqref{eq: V4 equations} there will be three additional symbol relations, which are
\begin{equation}\label{eq: secondary integrability symbol relations}
    \ftablB{2}{3}{2} = 0
    \qquad \ftablB{2}{3'}{2} = 0
    \qquad \ftablC{2}{3'}{0} - 2\ftablB{1}{3}{2} + 4\ftablB{2}{3}{1} = 0 .
\end{equation}
These are again most simply seen by the embedding into tableau forms described in the previous section.

\subsection{The additional integrability conditions}\label{sec:The additional integrability conditions}
Having made the reductions of Proposition \ref{thm: prop normalization of tableau}, a further integrability condition becomes apparent.
Indeed, one computes that on \( \V_{1} \),
\[ \d\sTheta{2}{2} \equiv 2\tablB{2}{3}{2} \theta^{2}\w\omega^{1'} \mod{\theta^{1}, \omega^{0}, \sTheta{a}{b}, \tThetaSingle{I}} , \]
so that necessarily
\[ \tablB{2}{3}{2} = 0 . \]
Let \( \V_{2} \) be the submanifold of \( \V_{1} \) on which this equation holds.
Note that a single additional symbol relation will hold on \( \V_{2} \).

There remain two further integrability conditions in terms of the jet linearization functions.
These are at higher order than the previous reductions, and can be seen after a prolongation of the differential system \( (\V_{2}, \I) \), but the exposition is simpler if we directly demonstrate them now.
This can be done by partially prolonging the system, because the higher integrability conditions will be apparent already on this partial prolongation.

To this end, it is not difficult to check that on \( \V_{2} \), 
\begin{align*}
  \d\sTheta{2}{2} & \equiv \tfrac{2}{3} \theta^{1} \w \eta^{2}_{1}
  & \mod{\sTheta{a}{b}, \tThetaSingle{I}} \\
  \d\tTheta{1}{3}
  & \equiv
  -\left(2 \omega^{1'} + \eta^{3}_{3'}\right) \w \omega^{2'}
  & \mod{\theta^{1}, \theta^{2}, \omega^{0}, \sTheta{a}{b}, \tThetaSingle{I}} \\
  \d \tTheta{2}{3}
  & \equiv \left(\tablB{1}{3}{2} \theta^{2} + \omega^{0}\right)\w\eta^{2}_{1} + \left(\tablB{2}{3'}{2} \theta^{2} + \omega^{1'}\right)\w\eta^{3}_{3'}
  & \mod{\sTheta{a}{b}, \tThetaSingle{I}} .
\end{align*}
From the first equation, the pullback of \( \eta^{2}_{1} \) by any integral section vanishes modulo the pullback of \( \theta^{1} \);
algebraically, we are free to calculate modulo \( \theta^{1} \w \eta^{2}_{1} \).
From this fact and the last two equations, the combination
\[ \eta^{3}_{3'} + 2\omega^{1'} + 2\tablB{2}{3'}{2}\theta^{2} \]
also pulls back to vanish modulo \( \theta^{1} \).
Since this is the case for every integral section, one may add the form
\[ F = \theta^{1}\w\left(\eta^{3}_{3'} + 2\omega^{1'} + 2\tablB{2}{3'}{2}\theta^{2}\right) \]
to the differential system \( (\V_{2}, \I) \) without losing solutions.

Inclusion of \( F \) has the immediate differential consequence,
\begin{align*}
  \d F \equiv \theta^{1}\w\left(
-2 \omega^{2'}\w\gamma^{1}_{2} 
- \omega^{1'}\w\left(2 \zeta_{2} 
+ 2 \zeta_{1} 
- 4 \eta^{3}_{3}\right)\right) \mod{F, \theta^{2}, \omega^{0}, \d\Theta^{22}, \sTheta{a}{b}, \tThetaSingle{I}} , 
\end{align*}
while some experimentation, guided by the symbol relations, leads to the following combination,
\begin{multline*}
  4 \d \tTheta{2}{3'} \theta^{1}\w\theta^{2} 
+ 18 \d \sTheta{1}{2}\w\theta^{1}\w\omega^{2'} 
- 12 \d \sTheta{1}{2}\w\theta^{2}\w\omega^{1'} 
- 3 \d \sTheta{1}{1}\w\theta^{1}\w\omega^{1'}
\\\equiv \theta^{1}\w\theta^{2}\w\left(24 \tablB{1}{3}{2} \omega^{1'}\w\omega^{2'} 
+ 2 \omega^{2'}\w\gamma^{1}_{2} 
+ \omega^{1'}\w\left(2 \zeta_{2} 
+ 2 \zeta_{1} 
- 4 \eta^{3}_{3}\right)\right) 
\\  \mod{F, \omega^{0}, \d\Theta^{22}, \sTheta{a}{b}, \tThetaSingle{I}} .
\end{multline*}
Both must vanish on integral sections, so \( \tablB{1}{3}{2} \) is seen to be an integrability condition.
Restrict to the submanifold \( \V_{3} \subset \V_{2} \) on which it vanishes.

After this, one further integrability condition becomes apparent, from
\begin{multline*}
 3 \d \sTheta{1}{2}\w\theta^{1}\w\omega^{2'} 
 + 2 \d \tTheta{1}{3}\w\theta^{1}\w\omega^{0} 
+ 4 \d \tTheta{2}{3}\w\theta^{2}\w\omega^{0} 
+ \d \tTheta{2}{3'}\w\theta^{1}\w\theta^{2}
\\ \equiv 
-10 \tablB{2}{3'}{2} \theta^{1}\w\theta^{2}\w\omega^{0}\w\omega^{2'} 
 \mod{F, \omega^{1'}, \d\Theta^{22}, \sTheta{a}{b}, \tThetaSingle{I}} .
\end{multline*}
Let \( \V_{4} \subset \V_{3} \) be the \( 35 \)-dimensional submanifold on which \( \tablB{2}{3'}{2} = 0 \).
The additional symbol relations from restricting to \( \V_{4} \) have already been described, in Equation \eqref{eq: secondary integrability symbol relations}.

\subsection{The remaining torsion absorption}\label{sec:The remaining torsion absorption}
After restriction to \( \V_{4} \), the remaining torsion of the linear Pfaffian system is quite simple.
Indeed, we have already that
\[ \d\sTheta{a}{b}
   \equiv{} - \ftablA{a}{b}{c} \w \theta^{c}
   \mod{\sTheta{d}{e},\tThetaSingle{I}} \]
and for \( I \neq \tendex{1}{3} \)
\[ \d\tThetaSingle{I}
  \equiv{} - \ftablBSingle{I}{c} \w \theta^{c} - \ftablCSingle{I}{i} \w \omega^{i}
    \mod{\sTheta{d}{e},\tThetaSingle{I}} , \]
so that the only remaining torsion is the last term in
\[ \d\tTheta{1}{3}
  \equiv{} - \ftablB{1}{3}{c} \w \theta^{c} - \ftablC{1}{3}{i} \w \omega^{i} - 2\omega^{1'}\w\omega^{2'}
    \mod{\sTheta{d}{e},\tThetaSingle{I}} . \]
This torsion is easily absorbed, by replacing \( \ftablC{1}{3}{2'} \) with
\[ \hat{\pi}^{\tendex{1}{3}}_{2'} = \ftablC{1}{3}{2'} + 2\omega^{1'} ={} \eta^{3}_{3'} + 2\omega^{1'} . \]
The symbol relations require one to make the same modification for \( \ftablC{2}{3}{1'} \), so that
\[ \hat{\pi}^{\tendex{2}{3}}_{1'} = \ftablC{2}{3}{1'} + 2\omega^{1'} ={} \eta^{3}_{3'} + 2\omega^{1'} . \]

\section{Higher order reductions}\label{sec:Higher order reductions}
\subsection{The failure of involutivity}\label{sec:The failure of involutivity}
Now that the torsion of the Pfaffian system on \( \V_{4} \) has been absorbed, one should ask whether the system is involutive.
Unfortunately, it is not.
This can be checked simply on the linearized tableau bundle of \( (\V_{4}, \I) \), which is involutive if and only if the differential system is.
Per the discussion in Subsection \ref{sec:The first integrability condition}, the tableau bundle is easily read off from the \( \frg_{N,-} \otimes \frg_{M,-}^{\vee} \)-valued tableau form, which on \( \V_{4} \) is
\begin{align*}
  \left( \begin{matrix}
      \ftablA{1}{1}{1} & \ftablA{1}{1}{2} & 0 & 0 & 0 \\
      \ftablA{1}{2}{1} & \ftablA{1}{2}{2} & 0 & 0 & 0 \\
      \ftablA{2}{2}{1} & 0 & 0 & 0 & 0 \\
      \ftablB{1}{3}{1} & \ftablB{1}{3}{2} & \ftablC{1}{3}{0} & 0 & \ftablC{1}{3}{2'} \\
      \ftablB{1}{3'}{1} & \ftablB{1}{3'}{2} & \ftablC{1}{3'}{0} & \tfrac{3}{2}\ftablA{1}{1}{1} & \tfrac{3}{2}\ftablA{1}{1}{2} - \ftablC{1}{3}{0} \\
      \ftablB{2}{3}{1} & 0 & \tfrac{3}{2}\ftablA{2}{2}{1} & \ftablC{1}{3}{2'} & 0 \\
      \ftablB{2}{3'}{1} & 0 & 2\ftablB{1}{3}{2} - 4\ftablB{2}{3}{1} & 3\ftablA{1}{2}{1} - \ftablC{1}{3}{0} & 3\ftablA{1}{2}{2} - \tfrac{3}{2}\ftablA{2}{2}{1}
    \end{matrix} \right) .
\end{align*}
Denote by \( A \) the linearized tableau at any point, noting that it is independent of which point is chosen.
The Cartan characters of \( A \) are, by inspection, given by
\[ (s_{1}, s_{2}, s_{3}, s_{4}, s_{5}) = (7, 4, 2, 1, 0) , \]
so that by Cartan's test for involutivity, \cite{Bryant.Chern.eaExteriorDifferentialSystems1991},
the prolongation \( A^{(1)} \) would need to have dimension \( 25 \) to be involutive,
whereas a calculation shows \( A^{(1)} \) to be only \( 18 \) dimensional.
To make this calculation, take the following as a basis of the tableau \( A \subset \frg_{N,-} \otimes \frg_{M,-}^{\vee} \),
\begin{align*}
& e_{\symdex{1}{1}} \otimes e^{1} + \tfrac{3}{2}e_{\tendex{1}{3'}} \otimes e^{1'},\quad
e_{\symdex{1}{1}} \otimes e^{2} + \tfrac{3}{2}e_{\tendex{1}{3'}} \otimes e^{2'},\quad
e_{\symdex{1}{2}} \otimes e^{1} + 3e_{\tendex{2}{3'}} \otimes e^{1'}, \\
& e_{\symdex{1}{2}} \otimes e^{2} + 3e_{\tendex{2}{3'}} \otimes e^{2'},\quad
e_{\symdex{2}{2}} \otimes e^{1} + \tfrac{3}{2}e_{\tendex{2}{3}} \otimes e^{0} - \tfrac{3}{2}e_{\tendex{2}{3'}} \otimes e^{2'},\quad
e_{\tendex{1}{3}} \otimes e^{1}, \\
& e_{\tendex{1}{3}} \otimes e^{2} + 2e_{\tendex{2}{3'}} \otimes e^{0},\quad
e_{\tendex{1}{3}} \otimes e^{0} - e_{\tendex{1}{3'}} \otimes e^{2'} - e_{\tendex{2}{3'}} \otimes e^{1'},\quad
e_{\tendex{1}{3}} \otimes e^{2'} + e_{\tendex{2}{3}} \otimes e^{1'}, \\
& e_{\tendex{1}{3'}} \otimes e^{1},\quad
e_{\tendex{1}{3'}} \otimes e^{2},\quad
e_{\tendex{1}{3'}} \otimes e^{0},\quad
e_{\tendex{2}{3}} \otimes e^{1} - 4e_{\tendex{2}{3'}} \otimes e^{0},\quad
e_{\tendex{2}{3'}} \otimes e^{1} .
\end{align*}
Then it is linear algebra from the definition of \( A^{(1)} \),
\[ A^{(1)} = \left(A \otimes \frg_{M,-}^{\vee}\right) \cap \left(\frg_{N,-} \otimes \Sym^{2} \frg_{M,-}^{\vee}\right) , \]
to compute the 18 element basis
\begin{equation}\label{eq: jet tableau prolongation}
  \begin{aligned}
    & \cblock{\Hlltltl:} e_{\symdex{1}{1}} \otimes e^{1}\circ e^{1} + 3e_{\tendex{1}{3'}} \otimes e^{1}\circ e^{1'},
    \cblock{\Hlltltz:} \qquad 2e_{\symdex{1}{1}} \otimes e^{1}\circ e^{2} + 3e_{\tendex{1}{3'}} \otimes e^{1}\circ e^{2'} + 3e_{\tendex{1}{3'}} \otimes e^{2}\circ e^{1'}, \\
    & \cblock{\Hlltztz:} e_{\symdex{1}{1}} \otimes e^{2}\circ e^{2} + 3e_{\tendex{1}{3'}} \otimes e^{2}\circ e^{2'},
    \cblock{\Hlztltl:} \qquad e_{\symdex{1}{2}} \otimes e^{1}\circ e^{1} + 6e_{\tendex{2}{3'}} \otimes e^{1}\circ e^{1'}, \\
    & \cblock{\Hlztltz:} e_{\symdex{1}{2}} \otimes e^{1}\circ e^{2} + 3e_{\tendex{2}{3'}} \otimes e^{1}\circ e^{2'} + 3e_{\tendex{1}{3}} \otimes e^{2}\circ e^{0} - 3e_{\tendex{1}{3'}} \otimes e^{2}\circ e^{2'} + 3e_{\tendex{2}{3'}} \otimes e^{0}\circ e^{0}, \\
    & \cblock{\Hzztltl:} e_{\symdex{2}{2}} \otimes e^{1}\circ e^{1} + 3e_{\tendex{2}{3}} \otimes e^{1}\circ e^{0} - 3e_{\tendex{2}{3'}} \otimes e^{1}\circ e^{2'} - 6e_{\tendex{2}{3'}} \otimes e^{0}\circ e^{0},
    \cblock{\Hlbtltl:} \qquad e_{\tendex{1}{3}} \otimes e^{1}\circ e^{1}, \\
    & \cblock{\Hlbtltz:} e_{\tendex{1}{3}} \otimes e^{1}\circ e^{2} + 2e_{\tendex{2}{3'}} \otimes e^{1}\circ e^{0},
    \cblock{\Hlbtlwo:} \qquad e_{\tendex{1}{3}} \otimes e^{1}\circ e^{0} - e_{\tendex{1}{3'}} \otimes e^{1}\circ e^{2'} - e_{\tendex{2}{3'}} \otimes e^{1}\circ e^{1'}, \\
    & \cblock{\Hlbtlwz:} e_{\tendex{1}{3}} \otimes e^{1}\circ e^{2'} + e_{\tendex{2}{3}} \otimes e^{1}\circ e^{1'} + 2e_{\tendex{1}{3}} \otimes e^{0}\circ e^{0} - 4e_{\tendex{1}{3'}} \otimes e^{0}\circ e^{2'} - 4e_{\tendex{2}{3'}} \otimes e^{0}\circ e^{1'}, \\
    & \cblock{\Hlbptltl:} e_{\tendex{1}{3'}} \otimes e^{1}\circ e^{1},
    \cblock{\Hlbptltz:} \qquad e_{\tendex{1}{3'}} \otimes e^{1}\circ e^{2},
    \cblock{\Hlbptztz:} \qquad e_{\tendex{1}{3'}} \otimes e^{2}\circ e^{2},
    \cblock{\Hlbptlwo:} \qquad e_{\tendex{1}{3'}} \otimes e^{1}\circ e^{0},
    \\
    & \cblock{\Hlbptzwo:} e_{\tendex{1}{3'}} \otimes e^{2}\circ e^{0},
    \cblock{\Hlbpwowo:} \qquad e_{\tendex{1}{3'}} \otimes e^{0}\circ e^{0},
    \cblock{\Hzbtltl:} \qquad e_{\tendex{2}{3}} \otimes e^{1}\circ e^{1} - 8e_{\tendex{2}{3'}} \otimes e^{1}\circ e^{0},
    \cblock{\Hzbptltl:} \qquad e_{\tendex{2}{3'}} \otimes e^{1}\circ e^{1} .
  \end{aligned}
\end{equation}

\subsection{Prolongation}\label{sec:Prolongation}
Having computed \( A^{(1)} \), the prolongation \( (\V_{4}^{(1)}, \I^{(1)}) \) of \( (\V_{4}, \I) \) is straightforward to describe.
The manifold \( \V_{4}^{(1)} \) is simply \( \V_{4} \times A^{(1)} \).
The basis of \( A^{(1)} \) given in display \eqref{eq: jet tableau prolongation} induces linear coordinates on \( A^{(1)} \), the \emph{tableau coefficients}
\begin{multline*}
  \Hlltltl, \Hlltltz, \Hlltztz, \Hlztltl, \Hlztltz, \Hzztltl, \Hlbtltl, \Hlbtltz, \Hlbtlwo, \\
  \Hlbtlwz, \Hlbptltl, \Hlbptltz, \Hlbptztz, \Hlbptlwo, \Hlbptzwo, \Hlbpwowo, \Hzbtltl, \Hzbptltl ,
\end{multline*}
with the corresponding functions on \( \V_{4}^{(1)} \) denoted by the same symbols.
Define on \( \V_{4}^{(1)} \) the following forms, in accordance with \( A^{(1)} \),
\[  \left( \begin{matrix}
      \Theta^{\symdex{1}{1}}_{\thindex{1}} & \Theta^{\symdex{1}{1}}_{\thindex{2}} \\
      \Theta^{\symdex{1}{2}}_{\thindex{1}} & \Theta^{\symdex{1}{2}}_{\thindex{2}} \\
      \Theta^{\symdex{2}{2}}_{\thindex{1}} & 0
    \end{matrix} \right)
  = \left( \begin{matrix}
      \pi^{\symdex{1}{1}}_{\thindex{1}} & \pi^{\symdex{1}{1}}_{\thindex{2}} \\
      \pi^{\symdex{1}{2}}_{\thindex{1}} & \pi^{\symdex{1}{2}}_{\thindex{2}} \\
      \pi^{\symdex{2}{2}}_{\thindex{1}} & 0
\end{matrix} \right)
  + \left(
    \begin{matrix}
      \Hlltltl\theta^{1} + \Hlltltz\theta^{2} & \Hlltltz\theta^{1} + \Hlltztz\theta^{2} \\
      \Hlztltl\theta^{1} + \Hlztltz\theta^{2} & \Hlztltz\theta^{1} \\
      \Hzztltl\theta^{1} & 0
    \end{matrix} \right), \]
\begin{align*}
  \Theta^{\tendex{1}{3 }}_{\thindex{1}} & = \pi^{\tendex{1}{3 }}_{\thindex{1}} + (  \Hlbtltl\theta^{1} + \Hlbtltz\theta^{2} + \Hlbtlwo\omega^{0} + \Hlbtlwz\omega^{2'})
  \\ \Theta^{\tendex{1}{3'}}_{\thindex{1}} & = \pi^{\tendex{1}{3'}}_{\thindex{1}} + ( \Hlbptltl\theta^{1} + \Hlbptltz\theta^{2} + \Hlbptlwo\omega^{0} + \tfrac{3}{2}\Hlltltl\omega^{1'} + \tfrac{3}{2}\Hlltltz\omega^{2'} - \Hlbtlwo\omega^{2'})
  \\ \Theta^{\tendex{2}{3 }}_{\thindex{1}} & = \pi^{\tendex{2}{3 }}_{\thindex{1}} + ( \Hzbtltl\theta^{1} + \tfrac{3}{2}\Hzztltl\omega^{0} + \Hlbtlwz\omega^{1'})
  \\ \Theta^{\tendex{2}{3'}}_{\thindex{1}} & = \pi^{\tendex{2}{3'}}_{\thindex{1}} + ( \Hzbptltl\theta^{1} + 2\Hlbtltz\omega^{0} - 4\Hzbtltl\omega^{0} + 3\Hlztltl\omega^{1'} - \Hlbtlwo\omega^{1'} + 3\Hlztltz\omega^{2'} - \tfrac{3}{2}\Hzztltl\omega^{2'})
  \\ \Theta^{\tendex{1}{3 }}_{\thindex{2}} & = \pi^{\tendex{1}{3 }}_{\thindex{2}} + ( \Hlbtltz\theta^{1} + 3\Hlztltz\omega^{0})
  \\ \Theta^{\tendex{1}{3'}}_{\thindex{2}} & = \pi^{\tendex{1}{3'}}_{\thindex{2}} + ( \Hlbptltz\theta^{1} + \Hlbptztz\theta^{2} + \Hlbptzwo\omega^{0} + \tfrac{3}{2}\Hlltltz\omega^{1'} + \tfrac{3}{2}\Hlltztz\omega^{2'} - 3\Hlztltz\omega^{2'}) ,
\end{align*}
\begin{align*}
  \Theta^{\tendex{1}{3}}_{\omindex{0}} & = \pi^{\tendex{1}{3}}_{\omindex{0}} + (\Hlbtlwo\theta^{1} + 3\Hlztltz\theta^{2} + 4\Hlbtlwz\omega^{0})
  \\ \Theta^{\tendex{1}{3}}_{\omindex{2'}} & = \pi^{\tendex{1}{3}}_{\omindex{2'}} + 2\omega^{1'} + \Hlbtlwz\theta^{1}
  \\ \Theta^{\tendex{1}{3'}}_{\omindex{0}} & = \pi^{\tendex{1}{3'}}_{\omindex{0}} + (\Hlbptlwo\theta^{1} + \Hlbptzwo\theta^{2} + \Hlbpwowo\omega^{0} - 4\Hlbtlwz\omega^{2'}).
\end{align*}
The prolongation ideal is the differential ideal
\begin{equation*}
  \I^{(1)} = \left\{ \sTheta{a}{b}, \tThetaSingle{I}, \ \Theta^{\symdex{1}{1}}_{\thindex{1}}, \Theta^{\symdex{1}{1}}_{\thindex{2}}, \Theta^{\symdex{1}{2}}_{\thindex{1}}, \Theta^{\symdex{1}{2}}_{\thindex{2}},  \Theta^{\symdex{2}{2}}_{\thindex{1}}, \Theta^{\tendex{1}{3 }}_{\thindex{1}}, \Theta^{\tendex{1}{3'}}_{\thindex{1}}, \Theta^{\tendex{2}{3 }}_{\thindex{1}}, \Theta^{\tendex{2}{3'}}_{\thindex{1}}, \Theta^{\tendex{1}{3 }}_{\thindex{2}}, \Theta^{\tendex{1}{3'}}_{\thindex{2}}, \Theta^{\tendex{1}{3}}_{\omindex{0}}, \Theta^{\tendex{1}{3}}_{\omindex{2'}}, \Theta^{\tendex{1}{3'}}_{\omindex{0}} \right\} .
\end{equation*}

On the prolongation \( \V^{(1)}_{4} \), many of the coframing forms are determined as linear combinations of \( \theta^{a}, \omega^{i} \) modulo \( \I^{(1)} \), a reflection of the coframe reductions made when normalizing, Proposition \ref{thm: prop normalization of tableau}.
Because of this, we can write the ideal \( \I^{(1)} \) in terms of a more geometrically adapted set of generators.
Indeed, comparing the equations in Subsection \ref{sec:The tableau forms} and the definition of prolongation forms just given, one finds that the following forms are each contained in \( \I^{(1)} \),
{\allowdisplaybreaks
\begin{align*}
 \tilde{\gamma}^{1}_{2} := {}& \gamma^{1}_{2} - (3 \Hlztltz \theta^{1} - 3 \Hzztltl \theta^{1})
 \\  \tilde{\gamma}^{0}_{2} := {}&  \gamma^{0}_{2} - (\Hlbtltz \theta^{1} - 2 \Hzbtltl \theta^{1} + 3 \Hlztltz \omega^{0} - 3 \Hzztltl \omega^{0} - 2 \Hlbtlwz \omega^{1'})
 \\  \tilde{\gamma} := {}&  \gamma - (-\Hlbptltz \theta^{1} + 2 \Hzbptltl \theta^{1} - \Hlbptztz \theta^{2} + 4 \Hlbtltz \omega^{0} - \Hlbptzwo \omega^{0} - 8 \Hzbtltl \omega^{0} - \tfrac{3}{2} \Hlltltz \omega^{1'}
 \\ & + 6 \Hlztltl \omega^{1'} - 2 \Hlbtlwo \omega^{1'} - \tfrac{3}{2} \Hlltztz \omega^{2'} + 9 \Hlztltz \omega^{2'} - 3 \Hzztltl \omega^{2'})
 \\  \tilde{\eta}^{1}_{1} := {}&  \eta^{1}_{1} - (-\tfrac{3}{4} \Hlltltz \theta^{1} - \tfrac{3}{4} \Hlltztz \theta^{2} + \tfrac{1}{2} \zeta_{2} + \tfrac{3}{2} \zeta_{1})
 \\  \tilde{\eta}^{1}_{2} := {}&  \eta^{1}_{2} - (-\tfrac{3}{2} \Hlltltl \theta^{1} - \tfrac{3}{2} \Hlltltz \theta^{2} + \gamma^{2}_{1})
 \\  \tilde{\eta}^{2}_{1} := {}&  \eta^{2}_{1} - (-\tfrac{3}{2} \Hzztltl \theta^{1})
 \\  \tilde{\eta}^{2}_{2} := {}&  \eta^{2}_{2} - (\tfrac{3}{4} \Hlltltz \theta^{1} - 3 \Hlztltl \theta^{1} + \tfrac{3}{4} \Hlltztz \theta^{2} - 3 \Hlztltz \theta^{2} + \tfrac{3}{2} \zeta_{2} + \tfrac{3}{2} \zeta_{1})
 \\  \tilde{\eta}^{3}_{3} := {}&  \eta^{3}_{3} - (\tfrac{3}{4} \Hlltltz \theta^{1} - \Hlbtlwo \theta^{1} + \tfrac{3}{4} \Hlltztz \theta^{2} - 3 \Hlztltz \theta^{2} - 4 \Hlbtlwz \omega^{0} + \tfrac{1}{2} \zeta_{2} + \tfrac{1}{2} \zeta_{1})
 \\  \tilde{\eta}^{3}_{3'} := {}&  \eta^{3}_{3'} - (-\Hlbtlwz \theta^{1} - 2 \omega^{1'})
 \\  \tilde{\eta}^{3'}_{3} := {}&  \eta^{3'}_{3} - (-\Hlbptlwo \theta^{1} - \Hlbptzwo \theta^{2} - \Hlbpwowo \omega^{0} + 4 \Hlbtlwz \omega^{2'} - 2 \gamma^{0}_{1})
 \\  \tilde{\eta}_{\tendex{1}{3}} := {}&  \eta_{\tendex{1}{3}} - (-3 \Hlbptltz \theta^{1} + 3 \Hzbptltl \theta^{1} - 3 \Hlbptztz \theta^{2} + 6 \Hlbtltz \omega^{0} - 3 \Hlbptzwo \omega^{0} - 12 \Hzbtltl \omega^{0}
 \\ & - \tfrac{9}{2} \Hlltltz \omega^{1'} + 9 \Hlztltl \omega^{1'} - 3 \Hlbtlwo \omega^{1'} - \tfrac{9}{2} \Hlltztz \omega^{2'} + 18 \Hlztltz \omega^{2'} - \tfrac{9}{2} \Hzztltl \omega^{2'})
 \\  \tilde{\eta}_{\tendex{1}{3'}} := {}&  \eta_{\tendex{1}{3'}} - (-3 \Hzbtltl \theta^{1} - \tfrac{9}{2} \Hzztltl \omega^{0} - 3 \Hlbtlwz \omega^{1'})
 \\  \tilde{\eta}_{\tendex{2}{3}} := {}&  \eta_{\tendex{2}{3}} - (-3 \Hlbptltl \theta^{1} - 3 \Hlbptltz \theta^{2} - 3 \Hlbptlwo \omega^{0} - \tfrac{9}{2} \Hlltltl \omega^{1'} - \tfrac{9}{2} \Hlltltz \omega^{2'} + 3 \Hlbtlwo \omega^{2'})
 \\  \tilde{\eta}_{\tendex{2}{3'}} := {}&  \eta_{\tendex{2}{3'}} - (-3 \Hlbtltl \theta^{1} - 3 \Hlbtltz \theta^{2} - 3 \Hlbtlwo \omega^{0} - 3 \Hlbtlwz \omega^{2'} + 3 \gamma^{0}_{1}) .
\end{align*}
}
It is not difficult to check that furthermore
\begin{equation}\label{eq: prolonged coframing forms}
  \I^{(1)} = \left\{ \sTheta{a}{b}, \tThetaSingle{I}, \ \tilde{\gamma}^{1}_{2}, \tilde{\gamma}^{0}_{2}, \tilde{\gamma}, \tilde{\eta}^{1}_{1}, \tilde{\eta}^{1}_{2}, \tilde{\eta}^{2}_{1}, \tilde{\eta}^{2}_{2}, \tilde{\eta}^{3}_{3}, \tilde{\eta}^{3}_{3'}, \tilde{\eta}^{3'}_{3}, \tilde{\eta}_{\tendex{1}{3}}, \tilde{\eta}_{\tendex{1}{3'}}, \tilde{\eta}_{\tendex{2}{3}}, \tilde{\eta}_{\tendex{2}{3'}} \right\} .
\end{equation}
This system of equations will be simplified in short order.

\subsection{Higher reductions}\label{sec:Higher reductions}
We now describe a series of reductions of \( (\V_{4}^{(1)}, \I^{(1)}) \) that can be determined from the prolonged coframing forms.
In doing so, all \( 18 \) of the tableau coefficients can be eliminated, resulting in a submanifold \( \V' \subset \V^{(1)}_{4} \) for which the projection \( \V' \to \V_{4} \) is an isomorphism.
It will furthermore be seen that the prolongation manifold of \( (\V', \I^{(1)}) \) is again isomorphic to \( \V' \), but with larger ideal \( \I^{(2)} \).
In the end, the tableau of \( \I^{(2)} \) is trivial, so that the remaining integrability conditions are easily determined, which will be taken up in the following subsection.

The reductions come in two types.
The first type is by integrability conditions, all of which are linear relations between the tableau coefficients and the curvature functions of \( \B_{M} \).
The second type uses the remaining freedoms of the coframing group to absorb components of the tableau coefficients.
For an example of the first type, one can check using the forms given in Equation \eqref{eq: prolonged coframing forms} that
\begin{align*}
  \d\tilde{\eta}^{3}_{3'}
   \equiv {}& \left(3\Hlltztz - 12\Hlztltz\right) \theta^{2}\w\omega^{1'} - 15\Hlbtlwz \omega^{0}\w\omega^{1'} \mod{\theta^{1}, \I^{(1)}_{1}} ,
\end{align*}
so it is an integrability condition that \( \Hlltztz = 4\Hlztltz \) and \( \Hlbtlwz = 0 \).
Here \( \I^{(1)}_{1} \) denotes the algebraic ideal generated by the \( 1 \)-forms in \( \I^{(1)} \).
Let \( \V_{5} \) be the submanifold on which these hold.
For an example of the second type, on \( \V_{5} \) it holds that
\begin{align*}
  \tilde{\gamma}^{1}_{2} & = \gamma^{1}_{2} - \tfrac{3}{4} u \theta^{1} ,
\end{align*}
where we have set \( u = \Hlltztz - 4\Hzztltl \), and so
\[ \d\tilde{\gamma}^{1}_{2} \equiv
  \left(\tfrac{3}{4}u \zeta_{2} + \tfrac{9}{4}u\zeta_{1} + \gamma_{2} - \tfrac{3}{4}\d u\right) \w \theta^{1}
  \mod{\theta^{2}, \omega^{0}, \omega^{1'}, \omega^{2'} , \I^{(1)}_{1}} . \]
From 
\begin{equation}\label{eq: gamma2 reduction}
   \d u \equiv u \zeta_{2} + 3u\zeta_{1} + \tfrac{4}{3}\gamma_{2} \mod{\theta^{a}, \omega^{i}, \I^{(1)}} 
\end{equation}
it follows that \( u \) varies under the action of \( G^{(1)} \) by scalings and the translation given by the exponential of \( X_{\gamma_{2}} \).
For this reason, every point of \( \V_{4}^{(1)} \) is in the \( G^{(1)} \) orbit of a point for which \( u = 0 \), and we may assume without losing solutions that \( \Hlltztz = 4\Hzztltl \).
Let \( \V_{6} \) be the submanifold of \( \V_{5} \) cut out by this condition.
Implicit in this argument is the fact that the exponential of \( X_{\gamma_{2}} \) is contained in \( G^{(1)} \), which holds because the reductions made in Proposition \ref{thm: prop normalization of tableau} were all of lower order.

With \( u = 0 \), the Equation \eqref{eq: gamma2 reduction} suggests the possibility of adding a linear combination of \( \gamma_{2} \) and the \( \theta^{a}, \omega^{i} \) to \( \I^{(1)} \) without losing solutions.
This is true, and could be done now, but there will be four other reductions of this kind, and we will treat these all at once at the end of this subsection.

We proceed in this manner, with the equations that augur reductions given in Table \ref{eq: array of why reductions} and the corresponding reductions collected in Table \ref{eq: table primary reductions}.
Note that the equations rely each in turn on the previous reductions.
Let \( \V' \subset \V_{6} \) be the submanifold cut out by all of the equations in Table \ref{eq: table primary reductions}.
\begin{table}
\[ \begin{array}{|l|l|}
\hline
\rule{0pt}{16pt}
1 & \d \tilde{\eta}^{3}_{3'} \equiv - 2 \Hlbptztz \theta^{1} \w \theta^{2} + (4 \Hlbtltz - 2 \Hlbptzwo - 14 \Hzbtltl) \theta^{1} \w \omega^{0}  \hfill \mod{\I^{(1)}_{1}}
\\[2pt] & \qquad\quad + (12 \Hlztltl - 8 \Hlbtlwo) \theta^{1} \w \omega^{1'}  \\[4pt]
\hline
\rule{0pt}{16pt}
2 & \d \tilde{\gamma} \equiv (- 3 u \zeta_{2} - \tfrac{9}{2} u \zeta_{1} + \gamma_{1} + \tfrac{3}{2} du) \w \omega^{1'} \hfill \mod{\theta^{1}, \theta^{2}, \omega^{0}, \omega^{2'}, \I^{(1)}_{1}} 
\\[2pt] & \qquad\quad \mbox{where} \quad u = \Hlltltz - 2 \Hlztltl  \\[4pt]
\hline
\rule{0pt}{16pt}
3 & \d \tilde{\eta}^{3'}_{3} \equiv 2 \Hlbpwowo \omega^{1'} \w \omega^{2'} \hfill \mod{\theta^{1}, \theta^{2}, \omega^{0}, \I^{(1)}_{1}}  \\[4pt]
\hline
\rule{0pt}{16pt}
4 & \d \tilde{\eta}^{2}_{2} \equiv - 3 \Hlbtltz \theta^{2} \w \omega^{1'} - \tfrac{3}{2} A_{3} \theta^{2} \w \omega^{1'} \hfill \mod{\theta^{1}, \omega^{0}, \omega^{2'}, \I^{(1)}_{1}}  \\[4pt]
\hline
\rule{0pt}{16pt}
5 & \tfrac{1}{14} \d \tilde{\eta}^{3}_{3} + \tfrac{1}{6}  \d \tilde{\eta}^{2}_{2} \equiv - (\Hzbtltl + \tfrac{2}{7} A_{3}) \theta^{2} \w \omega^{1'} \hfill \mod{\theta^{1}, \omega^{0}, \omega^{2'}, \I^{(1)}_{1}}  \\[4pt]
\hline
\rule{0pt}{16pt}
6 & 6 \d \tilde{\gamma}^{0}_{2} \w \theta^{2} - 4 \d \tilde{\eta}_{\tendex{1}{3'}} \w \theta^{2} + 3 \d \tilde{\eta}^{3'}_{3} \w \theta^{1} \equiv 30 \Hlbptltz \theta^{1} \w \theta^{2} \w \omega^{1'} \hfill \mod{\omega^{0}, \omega^{2'}, \I^{(1)}_{1}}  \\[4pt]
\hline
\rule{0pt}{16pt}
7 & \tfrac{2}{7} \d \tilde{\gamma}^{0}_{2} - \tfrac{1}{42} \d \tilde{\eta}_{\tendex{1}{3'}}  \equiv - (\Hzbptltl + \tfrac{2}{7} B_{3}) \theta^{1} \w \omega^{1'} \hfill \mod{\theta^{2}, \omega^{0}, \omega^{2'}, \I^{(1)}_{1}}  \\[4pt]
\hline
\rule{0pt}{16pt}
8 & \d \tilde{\eta}^{1}_{1} - \d \tilde{\eta}^{2}_{2} + 2 \d \tilde{\eta}^{3}_{3} \equiv (9 \Hlbtltl + 4 \Hlbptlwo + A_{4}) \theta^{1} \w \omega^{1'} \hfill \mod{\theta^{2}, \I^{(1)}_{1}}
\\[2pt] & \d \tilde{\eta}^{3'}_{3} \w \theta^{1} + 4 (\d \tilde{\eta}^{1}_{1} - \d \tilde{\eta}^{2}_{2} + \d \tilde{\eta}^{3}_{3}) \w \omega^{0} \equiv (- 24 \Hlbtltl + \Hlbptlwo - 6 A_{4}) \theta^{1} \w \omega^{0} \w \omega^{1'}  \\[4pt]
\hline
\rule{0pt}{16pt}
9 & \tfrac{5}{42} \d \tilde{\eta}^{3'}_{3} + \tfrac{2}{42} \d \tilde{\eta}_{\tendex{2}{3'}} \equiv (\Hlbptltl - \tfrac{2}{21} B_{4}) \theta^{1} \w \omega^{1'} \hfill \mod{\theta^{2}, \omega^{0}, \omega^{2'}, \I^{(1)}_{1}}  \\[4pt]
\hline
\rule{0pt}{16pt}
10 & \d \tilde{\eta}^{2}_{1} \equiv (-\tfrac{3}{2} \Hzztltl \zeta_{2} - \tfrac{9}{2} \Hzztltl \zeta_{1} + \tfrac{1}{3} \eta_{11} + \tfrac{3}{2} \d\Hzztltl) \w \theta^{1} \hfill \mod{\theta^{2}, \omega^{i}, \I^{(1)}_{1}}  \\[4pt]
\hline
\rule{0pt}{16pt}
11 & 3\d \tilde{\eta}^{1}_{1} - \d \tilde{\eta}^{2}_{2} \equiv (-6 \Hlztltl \zeta_{2} - 9 \Hlztltl \zeta_{1} + \tfrac{2}{3} \eta_{12} + 3 \d\Hlztltl) \w \theta^{1} \hfill \mod{\theta^{2}, \omega^{i}, \I^{(1)}_{1}}  \\[4pt]
\hline
\rule{0pt}{16pt}
12 & \d \tilde{\eta}^{1}_{2} \equiv (-\tfrac{9}{2} \Hlltltl \zeta_{2} - \tfrac{9}{2} \Hlltltl \zeta_{1} + \tfrac{1}{3} \eta_{22} + \tfrac{3}{2} \d\Hlltltl) \w \theta^{1} \hfill \mod{\theta^{2}, \omega^{i}, \I^{(1)}_{1}}  \\[4pt]
\hline
\end{array} \]
\caption{}
\label{eq: array of why reductions}
\[ \begin{array}{|l|l|l||l|l|l|}
\hline
 \rule{0pt}{16pt} 1 & \Hlbptztz = 0, \enspace \Hlbtlwo = \tfrac{3}{2}\Hlztltl &
 & 7 & \Hzbptltl = -\tfrac{2}{7} B_{3} & \\[2pt]
 & \enspace \Hlbptzwo = 2\Hlbtltz - 7\Hzbtltl & & & & \\[4pt]
 \hline
  \rule{0pt}{16pt} 2 & \Hlltltz = 2\Hlztltl & \gamma_{1}
  & 8 & \Hlbtltl = -\tfrac{5}{21} A_{4}, \enspace \Hlbptlwo = \tfrac{2}{7} A_{4} & \\[4pt]
  \hline
  \rule{0pt}{16pt} 3 & \Hlbpwowo = 0 &
  & 9 & \Hlbptltl = \tfrac{2}{21} B_{4} & \\[4pt]
  \hline
  \rule{0pt}{16pt} 4 & \Hlbtltz = -\tfrac{1}{2} A_{3} &
  & 10 & \Hlltztz = 0 & \eta_{11} \\[4pt]
  \hline
  \rule{0pt}{16pt} 5 & \Hzbtltl = -\tfrac{2}{7} A_{3} &
  & 11 & \Hlztltl = 0 & \eta_{12} \\[4pt]
  \hline
  \rule{0pt}{16pt} 6 & \Hlbptltz = 0 &
  & 12 & \Hlltltl = 0 & \eta_{22} \\[4pt]
  \hline
 \end{array} \]
 \caption{}
 \label{eq: table primary reductions}
\end{table}


Reducing to \( \V' \) involves 5 coframe reductions, after Equation \eqref{eq: gamma2 reduction} and in rows \( 2, 10, 11, 12 \) of Table \ref{eq: table primary reductions}.
Let \( G^{(2)} \) denote the subgroup of \( G^{(1)} \) that preserves these reductions, which can at this point be seen to be isomorphic to the \( 4 \)-dimensional subgroup \( G^{(2)}_{M} \) of \( G^{0}_{M} \).
Corresponding to the coframe reductions, there are \( 5 \) more \( 1 \)-forms that can be added to \( \I^{(1)} \) at no cost---they will automatically vanish for any solution.
Formally, these forms are of higher order, living on the prolongation \( \V^{(2)} \) of \( \V' \).
However, the prolongation manifold \( \V^{(2)} \) is isomorphic to \( \V' \), so it is only the prolongation ideal \( \I^{(2)} \) that increases.

Preliminary to determining these forms, note that on \( \V' \) it holds
\begin{equation}\label{eq: ABC curvature reductions}
\begin{aligned}
\d \tilde{\gamma}^{1}_{2} \w \omega^{0}  
+ \d \tilde{\gamma}^{0}_{2} \w \theta^{1}
\equiv {} 3 C_{1} \theta^{1} \w \theta^{2} \w \omega^{0} 
- B_{2} \theta^{1} \w \theta^{2} \w \omega^{1'} 
\hfill \mod{\I^{(1)}_{1}}
\\ 
- A_{2} \theta^{2} \w \omega^{0} \w \omega^{1'} 
- B_{1} \theta^{1} \w \theta^{2} \w \omega^{2'} 
- A_{2} \theta^{1} \w \omega^{0} \w \omega^{2'} 
- A_{1} \theta^{2} \w \omega^{0} \w \omega^{2'} ,
\end{aligned}
\end{equation}
so that each coefficient determines an integrability condition, and thus it is necessary to restrict to the submanifold of \( \V' \) where \( A_{1} = A_{2} = B_{1} = B_{2} = C_{1} = 0 \).
This is of course an empty set for the generic \( (2,3,5) \)-manifold, so we find the first curvature restrictions on embeddings.
There are several more restrictions along these lines, which will be taken up systematically in the following subsection.
For now, simply suppose that these \( 5 \) curvature functions vanish.

Returning to the prolongation forms, it now holds that
\begin{align*}
\d \tilde{\gamma}^{1}_{2} = {}& \theta^{1} \w \left(\tfrac{17}{14} A_{3} \omega^{1'} 
- \gamma_{2}\right) & \mod{\I^{\left(1\right)}_{1}}
\\ \d \tilde{\gamma}^{0}_{2} \equiv {}& \omega^{0} \w \left(
-2 C_{2} \theta^{1} 
- \tfrac{1}{14} A_{3;0} \theta^{1} 
- \gamma_{2}\right) & \mod{\theta^{2}, \omega^{1'} ,\omega^{2'}, \I^{\left(1\right)}_{1}} ,
\end{align*}
so \( \gamma_{2} \) is uniquely determined on solutions, and we may add
\[ \tilde{\gamma}_{2} = \gamma_{2}
+ 2 C_{2} \theta^{1} 
+ \tfrac{1}{14} A_{3;0} \theta^{1} 
- \tfrac{17}{14} A_{3} \omega^{1'} \]
to \( \I^{(1)} \).
There are four more reductions in this manner.
The equations in each row of Table \ref{eq: table prolongation reasons} determine uniquely one of \( \tilde{\gamma}_{1}, \tilde{\eta}_{11}, \tilde{\eta}_{12}, \tilde{\eta}_{22} \),
which is then given in Table \ref{eq: table prolongation results}.
The prolongation ideal \( \I^{(2)} \) is generated by these forms and \( \I^{(1)} \).
\begin{table}
\[ \begin{array}{|l|l|}
\hline
1 & 2 \d \tilde{\eta}^{3}_{3} \equiv \theta^{1} \w (
-C_{2} \theta^{2} 
+ \tfrac{22}{7} B_{3} \omega^{0} 
- \tfrac{9}{7} A_{4} \omega^{1'} 
- \tfrac{37}{14} A_{3} \omega^{2'} 
+ \gamma_{1}) \hfill \mod{\I^{(1)}_{1}} 
 \\ & \d \tilde{\gamma} \equiv \omega^{1'} \w (C_{3} \theta^{1} 
+ \tfrac{4}{7} B_{3;1'} \theta^{1} 
- \gamma_{1}) \hfill \mod{\theta^{2}, \omega^{0} ,\omega^{2'}, \I^{(1)}_{1}}  \\[4pt]
\hline
\rule{0pt}{16pt}
2 & 3 \d \tilde{\eta}^{2}_{1} \equiv \theta^{1} \w (\tfrac{18}{7} A_{3} \omega^{1'} 
- \eta_{11}) \hfill \mod{\I^{(1)}_{1}} 
 \\ & \d \tilde{\eta}_{\tendex{1}{3'}} \equiv \omega^{0} \w (
-\tfrac{6}{7} A_{3;0} \theta^{1} 
- \eta_{11}) \hfill \mod{\theta^{2}, \omega^{1'} ,\omega^{2'}, \I^{(1)}_{1}}  \\[4pt]
\hline
\rule{0pt}{16pt}
3 & 3 \d \tilde{\eta}^{2}_{2} 
- 9\d \tilde{\eta}^{3}_{3} \equiv \theta^{1} \w (
-\tfrac{54}{7} B_{3} \omega^{0} 
+ \tfrac{24}{7} A_{4} \omega^{1'} 
+ \tfrac{54}{7} A_{3} \omega^{2'} 
- \eta_{12}) \hfill \mod{\I^{(1)}_{1}} 
 \\ & \d \tilde{\eta}_{\tendex{1}{3}} \equiv \omega^{1'} \w (\tfrac{6}{7} B_{3;1'} \theta^{1} 
- \eta_{12}) \hfill \mod{\theta^{2}, \omega^{0} ,\omega^{2'}, \I^{(1)}_{1}}  \\[4pt]
\hline
\rule{0pt}{16pt}
4 & 3 \d \tilde{\eta}^{1}_{2} \equiv \theta^{1} \w (
-\tfrac{36}{7} B_{4} \omega^{0} 
+ 3 A_{5} \omega^{1'} 
+ \tfrac{36}{7} A_{4} \omega^{2'} 
- \eta_{22}) \hfill \mod{\I^{(1)}_{1}} 
 \\ & \d \tilde{\eta}_{\tendex{2}{3}} \equiv \omega^{1'} \w (\tfrac{2}{7} B_{4;1'} \theta^{1} 
- \eta_{22}) \hfill \mod{\theta^{2}, \omega^{0} ,\omega^{2'}, \I^{(1)}_{1}}  \\[4pt]
\hline
\end{array} \]
  \caption{}
  \label{eq: table prolongation reasons}
\[  \begin{array}{|l|l|}
  \hline
\rule{0pt}{16pt}
1 & \tilde{\gamma}_{1} = \gamma_{1} 
- C_{3} \theta^{1} 
- \tfrac{4}{7} B_{3;1'} \theta^{1} 
- C_{2} \theta^{2} 
+ \tfrac{22}{7} B_{3} \omega^{0} 
- \tfrac{9}{7} A_{4} \omega^{1'} 
- \tfrac{37}{14} A_{3} \omega^{2'}   \\[4pt]
\hline
\rule{0pt}{16pt}
2 & \tilde{\eta}_{11} = \eta_{11} 
+ \tfrac{6}{7} A_{3;0} \theta^{1} 
- \tfrac{18}{7} A_{3} \omega^{1'}  \\[4pt]
\hline
\rule{0pt}{16pt}
3 & \tilde{\eta}_{12} = \eta_{12} 
- \tfrac{6}{7} B_{3;1'} \theta^{1} 
+ \tfrac{54}{7} B_{3} \omega^{0} 
- \tfrac{24}{7} A_{4} \omega^{1'} 
- \tfrac{54}{7} A_{3} \omega^{2'}  \\[4pt]
\hline
\rule{0pt}{16pt}
4 & \tilde{\eta}_{22} = \eta_{22} 
-\tfrac{2}{7} B_{4;1'} \theta^{1} 
+ \tfrac{36}{7} B_{4} \omega^{0} 
- 3 A_{5} \omega^{1'} 
- \tfrac{36}{7} A_{4} \omega^{2'}
\\[4pt]
\hline
  \end{array} \]
    \caption{}
    \label{eq: table prolongation results}
\end{table}

The consequences of all the reductions made in this subsection are summarized in the statement that the following forms generate \( \I^{(2)} \):
\[ \sTheta{a}{b}, \tThetaSingle{I}, \]
\begin{equation}\label{eq: BP full reduction M}
  \begin{aligned}
    \tilde{\gamma}^{1}_{2} & = \gamma^{1}_{2},
    \qquad \tilde{\gamma}^{0}_{2} = \gamma^{0}_{2} - \tfrac{1}{14} A_{3} \theta^{1},
    \qquad \tilde{\gamma} = \gamma + \tfrac{4}{7} B_{3} \theta^{1} + \tfrac{5}{7} A_{3} \omega^{0},
    \\  \tilde{\gamma}_{2} & = \gamma_{2} + 2 C_{2} \theta^{1} + \tfrac{1}{14} A_{3;0} \theta^{1} - \tfrac{17}{14} A_{3} \omega^{1'},
    \\  \tilde{\gamma}_{1} & = \gamma_{1} - C_{3} \theta^{1} - \tfrac{4}{7} B_{3;1'} \theta^{1} - C_{2} \theta^{2} + \tfrac{22}{7} B_{3} \omega^{0} - \tfrac{9}{7} A_{4} \omega^{1'} - \tfrac{37}{14} A_{3} \omega^{2'},
  \end{aligned}
\end{equation}
\begin{equation}\label{eq: BP full reduction N}
  \begin{aligned}
    \left(
      \begin{matrix}
        \tilde{\eta}^{1}_{1} & \tilde{\eta}^{1}_{2} \\
        \tilde{\eta}^{2}_{1} & \tilde{\eta}^{2}_{2}
      \end{matrix} \right)
    & =
    \left(
      \begin{matrix}
        \eta^{1}_{1} & \eta^{1}_{2} \\
        \eta^{2}_{1} & \eta^{2}_{2}
      \end{matrix} \right)
    - \left(
      \begin{matrix}
        \tfrac{1}{2} \zeta_{2} + \tfrac{3}{2} \zeta_{1} &
        \gamma^{2}_{1} \\
        0 &
        \tfrac{3}{2} \zeta_{2} + \tfrac{3}{2} \zeta_{1}
      \end{matrix} \right)
    \\[3pt]
    \left(
      \begin{matrix}
        \tilde{\eta}^{3}_{3} & \tilde{\eta}^{3'}_{3} \\
        \tilde{\eta}^{3}_{3'} & - \tilde{\eta}^{3}_{3}
      \end{matrix} \right)
    & =
    \left(
      \begin{matrix}
        \eta^{3}_{3} & \eta^{3'}_{3} \\
        \eta^{3}_{3'} & - \eta^{3}_{3}
      \end{matrix} \right)
      - \left(
      \begin{matrix}
        \tfrac{1}{2} \zeta_{2} + \tfrac{1}{2} \zeta_{1} &
        -\tfrac{2}{7} A_{4} \theta^{1} - A_{3} \theta^{2} - 2 \gamma^{0}_{1} \\
        -2 \omega^{1'} &
        -\tfrac{1}{2} \zeta_{2} - \tfrac{1}{2} \zeta_{1}
      \end{matrix} \right)
    \\[3pt]
    \left(
      \begin{matrix}
        \tilde{\eta}_{\tendex{1}{3}} & \tilde{\eta}_{\tendex{1}{3'}} \\
        \tilde{\eta}_{\tendex{2}{3}} & \tilde{\eta}_{\tendex{2}{3'}}
      \end{matrix} \right)
    & = \left(
      \begin{matrix}
        \eta_{\tendex{1}{3}} & \eta_{\tendex{1}{3'}} \\
        \eta_{\tendex{2}{3}} & \eta_{\tendex{2}{3'}}
      \end{matrix} \right)
      - \left(
      \begin{matrix}
        -\tfrac{6}{7} B_{3} \theta^{1} - \tfrac{18}{7} A_{3} \omega^{0} &
        \tfrac{6}{7} A_{3} \theta^{1} \\
        -\tfrac{2}{7} B_{4} \theta^{1} - \tfrac{6}{7} A_{4} \omega^{0} &
        \tfrac{5}{7} A_{4} \theta^{1} + \tfrac{3}{2} A_{3} \theta^{2} + 3 \gamma^{0}_{1}
      \end{matrix} \right)
    \\[3pt]
    \left(
      \begin{matrix}
        \tilde{\eta}_{11} \\ \tilde{\eta}_{12} \\ \tilde{\eta}_{22}
      \end{matrix} \right)
    & = \left(
      \begin{matrix}
        \eta_{11} \\ \eta_{12} \\ \eta_{22}
      \end{matrix} \right)
    - \left(
      \begin{matrix}
        -\tfrac{6}{7} A_{3;0} \theta^{1} + \tfrac{18}{7} A_{3} \omega^{1'} \\
        \tfrac{6}{7} B_{3;1'} \theta^{1} - \tfrac{54}{7} B_{3} \omega^{0} + \tfrac{24}{7} A_{4} \omega^{1'} + \tfrac{54}{7} A_{3} \omega^{2'} \\
        \tfrac{2}{7} B_{4;1'} \theta^{1} - \tfrac{36}{7} B_{4} \omega^{0} + 3 A_{5} \omega^{1'} + \tfrac{36}{7} A_{4} \omega^{2'}
      \end{matrix} \right) .
  \end{aligned}
\end{equation}

\subsection{The curvature obstructions to embedding}\label{sec:The curvature obstructions to embedding}
After the reductions of the previous subsection, the linearized tableau of \( (\V', \I^{(2)}) \) is trivial. This can be seen by noting that \( \V' \) is \( 35 \)-dimensional and the independence condition generated by \( \I^{(2)} \) and \( \theta^{a}, \omega^{i} \) is of rank \( 31 \), while \( \I^{(2)} \) has \( 4 \) dimensions worth of Cauchy characteristic directions.
Because its tableau is trivial, \( \I^{(2)} \) is either a Frobenius ideal or it admits no integral sections.
In the former case, the existence of Pfaffian embeddings is clear, as they correspond to the leaves of the foliation induced by \( \I^{(2)} \), while in the latter case there are no embeddings.

Calculation of the obstructions to \( \I^{(2)} \) being Frobenius is long but routine.
The obstructions are entirely determined by the curvature of the \( (2,3,5) \)-manifold \( M \), and one finds about a hundred relations between the curvature functions and their semibasic derivatives.
However, almost all of the conditions hold if there exist coframing reductions to \( \B_{M} \) as stated in Theorem \ref{thm: main embeddability theorem}, whose proof we turn to now.

Consider the following diagram.
\[ \begin{tikzcd}[row sep={10mm,between origins}]
   (\V', \I^{(2)}) \ar[dr, "\cong" swap] \ar[r, hook] & \V^{(1)}_{4} \ar[d] & & \\
   & \V_{4} \ar[r, hook]\ar[d] & \BJets \ar[r] & \B_{M} \ar[dll] \\
   & M \ar[uul, bend left, dashed, "s_{1}"] \ar[urr, bend right, dashed, "s_{2}"] & & 
\end{tikzcd} \]
If an integral section \( s_{1} \) exists, then it pushes down to a section \( s_{2} \) of \( \B_{M} \), a complete reduction of \( \B_{M} \).
It holds by construction that \( s^{*}_{1} \Omega_{M} = s^{*}_{2} \Omega_{M} \), but it also holds that each form in \( \I^{(2)} \) pulls back to zero, so the reduction of \( \B_{M} \) satisfies in particular the conditions implied by Equations \eqref{eq: BP full reduction M}.

In case \( \B_{M} \) admits such a reduction, several relations on the curvature functions must hold, and in particular the following two.
First, on \( M \) with the coframing forms of \( \B_{M} \) pulled back by \( s_{2} \),
\[ 0 = \d \tilde{\gamma}_{2} \equiv (6 C_{2} - A_{3;0}) \omega^{0} \w \omega^{1'} \mod{\theta^{1}, \theta^{2}, \omega^{2'}} , \]
so that \( A_{3;0} = 6C_{2} \).
Second,
\[ 0 = 7 \d \tilde{\gamma}_{1} \w \theta^{2} + 9 \d^{2} (\gamma^{2}_{1}) \equiv -16(3 C_{3} + B_{3;1'}) \theta^{2} \w \omega^{0} \w \omega^{1'} \mod{\theta^{1}, \omega^{2'}} , \]
where we recall the final remark of Subsection \ref{sec:Curved (2,3,5)-manifolds}, and in particular Equation \eqref{eq: higher curvature relation example}.
So \( B_{3;1'} = -3 C_{3} \).
Substitution of these two relations into \eqref{eq: BP full reduction M} gives the reduction statement of Theorem \ref{thm: main embeddability theorem}.

One can compute all of the relations that follow from such a coframing reduction in a similar manner, and there are about a hundred all told.
The first order conditions have already been listed in Equation \eqref{eq: curvature consequences of reduction}.
In the end, one is obliged to restrict to the submanifold \( \V \subseteq \V' \) where these conditions on the curvature hold.
This can be thought of as a restriction to the \( 1 \)-jets whose adapted coframings are admissible as here.

Aside from the consequences of a coframing reduction, there are only two further conditions on \( \V \) for \( \I^{(2)} \) to be Frobenius.
The first is because of the equation
\[ \d \tilde{\eta}^{3'}_{3} \equiv 
-\left(\tfrac{10}{7} B_{4} + \tfrac{2}{7} A_{4;1'}\right) \theta^{1}\w\omega^{1'} \mod{\I^{(2)}} , \]
so that embeddability requires
\[ A_{4;1'} = -5 B_{4} . \]

The final integrability condition, from the differential of \( \tilde{\eta}_{22} \), is only seen at high enough order that special care is needed.
One finds that
\begin{equation}\label{eq: last integrability condition}
7 \d\tilde{\eta}_{22} \equiv
\begin{multlined}[t]
\theta^{1} \w \bigl(
-\tfrac{144}{7} A_{4} B_{3} \theta^{2} 
- \tfrac{60}{7} A_{3} B_{4} \theta^{2} 
+ \tfrac{60}{7} (A_{4})^2 \omega^{0} 
+ \tfrac{108}{7} A_{3} A_{5} \omega^{0} 
+ 36 B_{4;1} \omega^{0} 
\\
 + 21 A_{5;1} \omega^{1'} 
+ 108 \tilde{D}_{4} \omega^{2'} 
+ 5 A_{5;0} \zeta_{2} 
+ 6 A_{5;0} \zeta_{1} 
+ 36 C_{3} \gamma^{2}_{1} 
+ 24 B_{4} \gamma^{0}_{1} 
- \d A_{5;0}\bigr)
\\  \mod{\I^{(2)}_{1}} ,
\end{multlined}
\end{equation}
so that it is necessary to describe \( \d A_{5;0} \).
In fact, the vertical derivatives of \( A_{5;0} \), and many of its semibasic derivatives, are already determined.
This can be seen by considering the compatibility conditions from \( \d^{2} B_{4} = 0 \), which leads to
\begin{equation*}
 2 \d^2 B_{4} \equiv
 \begin{multlined}[t]
 \omega^{1'} \w \bigl(\tfrac{144}{7} A_{4} B_{3} \theta^{2} 
+ \tfrac{60}{7} A_{3} B_{4} \theta^{2} 
- \tfrac{60}{7} (A_{4})^2 \omega^{0} 
- \tfrac{108}{7} A_{3} A_{5} \omega^{0} 
- 36 B_{4;1} \omega^{0}
\\ 
- 108 \tilde{D}_{4} \omega^{2'} 
- 5 A_{5;0} \zeta_{2} 
- 6 A_{5;0} \zeta_{1} 
- 36 C_{3} \gamma^{2}_{1} 
- 24 B_{4} \gamma^{0}_{1} 
+ \d A_{5;0}\bigr)
\\ \mod{\theta^{1}, \I^{(2)}_{1}} ,
\end{multlined}
\end{equation*}
so that the only undetermined semibasic derivatives of \( A_{5;0} \) are \( A_{5;0;1} \) and \( A_{5;0;1'} \).
Comparison with Equation \eqref{eq: last integrability condition} results in the much simplified equation
\[ 7 \d\tilde{\eta}_{22} \equiv
\theta^{1} \w (21 A_{5;1} \omega^{1'} - A_{5;0;1'}\omega^{1'})
\qquad \mod{\I^{(2)}} . \]
If \( M \) satisfies \( A_{5;0;1'} = 21 A_{5;1} \), then \( \I^{(2)} \) is seen to be Frobenius, as desired.

To finish the proof of Theorem \ref{thm: main embeddability theorem}, note that if \( \I^{(2)} \) is Frobenius, then any maximal leaf projects to \( \B_{M} \) to give a coframe reduction, as required.
In fact, the leaf will be a \( G^{(2)}_{M} \) bundle over \( M \), and the reduction a \( G^{(2)}_{M} \)-principal reduction.

\section{Examples}\label{sec:Examples}
We finally give two examples of Pfaffian embeddings.
In both cases, the embedding is found using Theorem \ref{thm: main embeddability theorem}, but once given, they are straightforward to check directly.

\subsection{The flat (2,3,5)-model}\label{sec:The flat (2,3,5)-model}
Following Tanaka \cite{TanakaDifferentialSystemsGraded1970}, a \emph{fundamental graded algebra} is a finite dimensional, negatively graded (hence nilpotent) Lie algebra generated by its \( (-1) \)-graded component.
Given a fundamental graded algebra \( \frm \), the flat model \( M(\frm) \) is given by the \emph{standard differential system} of \( \frm \), the simply connected Lie group corresponding to \( \frm \), with distribution spanned by \( \frm_{-1} \).
Given fundamental graded algebras \( \frm, \frn \), symbol compatibility, meaning the existence of an injective graded Lie algebra morphism \( f \colon \frm \to \frn \), implies the existence of a Pfaffian embedding between the flat models,
because the corresponding group embedding \( F \colon M(\frm) \to M(\frn) \) is by construction a Pfaffian embedding.

The flat \( (2,3,5) \)-manifold is \emph{not} flat when regarded as a \( (3,5) \)-manifold, so this construction cannot quite determine a Pfaffian embedding to \( N \).
Nonetheless, a similar construction can be made.
The point is that while there cannot exist an injective filtered Lie algebra map of the symbol algebra \( \frg_{M,-} \) into the symbol algebra \( \frg_{N,-} \), there does exist an injective filtered Lie map into the prolongation \( \frg_{N} \).
Inspection of \( \I^{(2)} \) (Equations \eqref{eq: first framed contact forms}, \eqref{eq: BP full reduction N}, and \eqref{eq: BP full reduction M}), with all curvature functions zero, implies the Lie algebra map \( f \colon \frg_{M,-} \to \frg_{N} \) defined on the basis of \( \frg_{M,-} \) by
\begin{align*}
   e_{1'} & \mapsto e_{\tendex{2}{3'}} - 2X_{\eta^{3}_{3'}}, &
   e_{0} & \mapsto e_{\tendex{1}{3}}, &
   e_{1} & \mapsto \tfrac{1}{3} e_{\symdex{1}{2}}, \\
   e_{2'} & \mapsto e_{\tendex{1}{3'}}, & &&
   e_{2} & \mapsto \tfrac{2}{3} e_{\symdex{1}{1}} .
\end{align*}
This is by inspection a filtration preserving Lie algebra map.
With this map, one can use the Lie exponential to define a Pfaffian embedding of the flat \( (2,3,5) \)-manifold \( M \) into \( N \), as in the following diagram.
\[ \begin{tikzcd}
  \frg_{M,-} \ar[rr, "f"]\ar[d, "\cong", "\exp" swap] & & \frg_{N} \ar[d, "\exp" swap] \\
  M \ar[r, dashed, "F"] & N & G_{N} \ar[l] 
\end{tikzcd} \]

\subsection{A non-flat example}\label{sec:A non-flat example}

If one supposes the existence of an embeddable \( (2,3,5) \)-manifold \( M \) whose Cartan quartic has two distinct double roots, then Theorem \ref{thm: main embeddability theorem} leads one to the following reduction of the Cartan geometry,
\begin{equation}\label{eq: D6 coframe reduction}
  \begin{aligned}
    \gamma^{1}_{2} & = 0
    & \qquad & \zeta_{1} = -\tfrac{3}{2} \omega^{0} - \tfrac{1}{2} \zeta_{2}
    \\ \gamma^{2}_{1} & = 0
    && \gamma^{0}_{2} = \tfrac{1}{2} \theta^{1}
    \\ \gamma^{0}_{1} & = 6 \omega^{2'} + \tfrac{17}{2} \theta^{2}
    && \gamma = -5 \omega^{0}
    \\ \gamma_{2} & = -17 \theta^{1} + \tfrac{17}{2} \omega^{1'}
    && \gamma_{1} = 7 \theta^{2} + \tfrac{37}{2} \omega^{2'} .
  \end{aligned}
\end{equation}
The additional reductions of \( \zeta_{1}, \gamma^{2}_{1} \) and \( \gamma^{0}_{1} \) beyond the statement of the theorem are differential consequences of the respective normalizations \( A_{3} = 7 \) and \( A_{4} = B_{3} = 0 \), which can be made under the assumption that \( A_{3} \neq 0 \), seen for example by comparing table (A.5) of \cite{TheCartantheoreticClassificationMultiplytransitive2022}.
After this normalization, \( A_{5} = 0 \) is automatic from the assumption that the Cartan quartic has two double roots.
With these normalizations there are enough compatibility conditions on the framing that the derivatives of all curvature functions collapse to zero.
In the end, the only nonzero curvatures are
\[ A_{3} = 7,\quad C_{2} = \pm 7,\quad E = \tilde{E}_{2} = \tfrac{63}{2} . \]
In particular, although homogeneity was not assumed, it follows a posteriori that the unique embeddable \( (2,3,5) \)-system with two double roots is homogeneous.
(The sign of \( C_{2} \) can be made positive up to isomorphism.)

Substitution of these reductions into the structure equations of \( \B_{M} \) (refer to Equations \eqref{eq: 235 first order structure} and \eqref{eq: 235 higher order structure}) results in the following differential system,
\begin{equation}\label{eq: D6 MC forms}
  \begin{aligned}
  \d \theta^{1} & = -\tfrac{9}{2} \theta^{1} \w \omega^{0} + 3 \omega^{0} \w \omega^{1'} + \tfrac{1}{2} \theta^{1} \w \zeta_{2}
  \\ \d \theta^{2} & = -\tfrac{9}{2} \theta^{2} \w \omega^{0} + 3 \omega^{0} \w \omega^{2'} - \tfrac{1}{2} \theta^{2} \w \zeta_{2}
  \\ \d \omega^{0} & = 6 \theta^{1} \w \omega^{2'} + 8 \theta^{1} \w \theta^{2} + 2 \omega^{1'} \w \omega^{2'}
  \\ \d \omega^{1'} & = \tfrac{3}{2} \omega^{0} \w \omega^{1'} - 6 \theta^{1} \w \omega^{0} + \tfrac{1}{2} \omega^{1'} \w \zeta_{2}
  \\ \d \omega^{2'} & = -\tfrac{21}{2} \omega^{0} \w \omega^{2'} + 12 \theta^{2} \w \omega^{0} - \tfrac{1}{2} \omega^{2'} \w \zeta_{2}
  \\ \d \zeta_{2} & = 24 \theta^{1} \w \theta^{2} - 18 \omega^{1'} \w \omega^{2'} + 48 \theta^{2} \w \omega^{1'} - 6 \theta^{1} \w \omega^{2'} .
  \end{aligned}
\end{equation}
Because the coefficients are constant, we can build in the usual manner a homogeneous \( (2,3,5) \)-manifold \( M \).
Indeed, the forms \eqref{eq: D6 MC forms} are easily checked to be differentially compatible, and this is equivalent to the statement that they can be integrated to the left invariant coframing
\begin{equation}\label{eq: D6 coframing}
  \theta^{1}, \theta^{2}, \omega^{0}, \omega^{1'}, \omega^{2'},\quad \zeta_{2}
\end{equation}
of a simply connected Lie group \( \hat{G} \).
Let \( M \) be the homogeneous manifold \( M = \hat{G} / L \), where \( L \) is the \( 1 \)-dimensional subgroup generated by the annihilator of \( \theta^{1}, \theta^{2}, \omega^{0}, \omega^{1'}, \omega^{2'} \).
The corank \( 3 \)-distribution on \( \hat{G} \) annihilating \( \{ \theta^{1}, \theta^{2}, \omega^{0} \} \) is \( L \) invariant, so descends to an invariant \( 2 \)-distribution on \( M \).

Once found, it is simple to check that \( M \) is embeddable into \( N \).
To do so, one simply has to check that the system \( \I^{(2)} \) (refer to the end of Subsection \ref{sec:Higher reductions}) is Frobenius after substitution of the reductions made here.
The integral leaves each provide the graph of an embedding \( \varphi \colon \hat{G} \hookrightarrow \B_{N} = G_{N} \), and are each equivalent under the action of \( G_{N} \).
Supposing without loss that \( \varphi \) sends the neutral element of \( \hat{G} \) to the neutral element of \( G_{N} \), then \( \varphi \) is a group homomorphism, which is furthermore seen to send \( L \) into \( G^{0}_{N} \).
As such, the quotient \( M = \hat{G} / L \) embeds into \( N = G_{N} / G^{0}_{N} \), and this embedding is by construction a Pfaffian embedding.


The family of multiply-transitive two-double-root \( (2,3,5) \)-structures was studied by Cartan \cite{CartanSystemesPfaffCinq1910}, who demonstrated that they are parameterized by a single constant.
In \cite{Willse2019--HomogeneousRealDistributionsIsotropy} Willse describes the realization of many homogeneous two-double-root \( (2,3,5) \)-systems over \( \R \) as the rolling distributions of constant curvature surfaces, in either Riemannian or Lorentz signature.
For these models, Cartan's parameter is essentially the ratio \( \lambda \) of curvature of the two surfaces.
In Willse's Table 8, the first family \( D.6^{2-}_{\lambda} \)
is \( \sl(\R^{2}) \oplus \sl(\R^{2}) \) with isotropy fiber spanned by \( H + \tilde{H} \) and invariant rank-\( 3 \) distribution spanned by \( H + \tilde{H}, X - \lambda \tilde{X} \) and \( Y - \tilde{Y} \).
For positive \( \lambda \), the members of this family are the rolling distributions from two constant curvature Lorentz surfaces (de Sitter spaces) of curvature ratio \( \lambda \).

With the change of basis
\begin{align*}
  X^{*} & = -\theta^{1} + \omega^{1'}, &
  Y^{*} & = 12\theta^{2} + 3\omega^{2'}, &
  H^{*} & = \tfrac{3}{4}\omega^{0} + \tfrac{1}{4}\zeta_{2}, \\
  \tilde{X}^{*} & = -8\theta^{1} - 4\omega^{1'}, &
  \tilde{Y}^{*} & = -3\theta^{2} - 3\omega^{2'}, &
  \tilde{H}^{*} & = -\tfrac{15}{4}\omega^{0} + \tfrac{1}{4}\zeta_{2}
\end{align*}
it is easy to see that the structure equations
\begin{align*}
  \d X^{*} & = 2 X^{*} \w H^{*} &
  \d Y^{*} & = - 2 Y^{*} \w H^{*} &
  \d H^{*} & = - X^{*} \w Y^{*} \\
  \d \tilde{X}^{*} & = 2 \tilde{X}^{*} \w \tilde{H}^{*} &
  \d \tilde{Y}^{*} & = - 2 \tilde{Y}^{*} \w \tilde{H}^{*} &
  \d \tilde{H}^{*} & = - \tilde{X}^{*} \w \tilde{Y}^{*}
\end{align*}
hold, and thus that \( \operatorname{Lie}(\hat{G}) \) is isomorphic to \( \sl(\R^{2}) \oplus \sl(\R^{2}) \).
Furthermore, the annihilator of Willse's distribution, which is spanned by \( H^{*} - \tilde{H}^{*}, \lambda X^{*} + \tilde{X}^{*}, Y^{*} + \tilde{Y}^{*} \), agrees with the span of
\( \omega^{0}, \theta^{1}, \theta^{2} \) only when \( \lambda = 4 \).
In this case, the isotropy is also seen to agree with \( H + \tilde{H} \).
In other words, the embeddable is given by the rolling distribution of two de Sitter spaces of curvature ratio exactly \( 4 \).

\bibliographystyle{amsplain}
\bibliography{refs}
\end{document}